\documentclass[12pt,epsfig]{article}
\usepackage{amsfonts, amssymb}
\usepackage{epsfig}
\usepackage{latexsym}
\def\Proof{{\bf Proof}\hspace{.3in}}
\def\eop{\hspace*{\fill}$\Box$ \vskip \baselineskip}
\newfont{\Bb}{msbm10}
\def\ni{\noindent}
\def\be{\begin{enumerate}}
\def\ee{\end{enumerate}}
\def\ec{\end{center}}
\def\bc{\begin{center}}

\def\-1{{^{-1}}}

\def\g{{\gamma}}

\def\l{{\lambda}}
\def\p{{\partial}}
\def\CC{{\mathbb{C}}}

\def\RR{{\mathbb{R}}}

\def\ZZ{\mbox{\Bb Z}}
\def\tr{\pitchfork}
\def\Der{\mbox{Der}}

\def\AA{{\cal A}}

\def\OO{{\cal O}}

\def\to{{\ \rightarrow\ }}

\def\dim{{\mbox{dim}}}
\newtheorem{teo}{Theorem}[section]
\newtheorem{example}[teo]{Example}
\newtheorem{lema}[teo]{Lemma}
\newtheorem{prop}[teo]{Proposition}
\newtheorem{corl}[teo]{Corollary}
\newtheorem{obs}[teo]{Remark}
\begin{document}

\title{Vanishing topology of codimension 1 multi-germs over $\Bbb{R}$ and
$\Bbb{C}$
\footnote{
{\em Secondary} 14P25, 32S70}
}
\author{T. Cooper, D. Mond and R. Wik Atique}

\maketitle
\begin{abstract} 
We construct all ${\cal A}_e$-codimension 1 multi-germs of maps
$(k^n,T)\to (k^p,0)$, with $n\geq p-1$, $(n,p)$ nice dimensions,
$k=\CC$ or $\RR$, by augmentation and concatenation operations, starting from
mono-germs ($|T|=1$) and one 0-dimensional bi-germ. 
As an application, we prove general statements for
multi-germs of corank $\leq 1$: every one has a real form with 
real perturbation carrying the vanishing homology of the complexification, 
every one is quasihomogeneous,
and when $n=p-1$ every one has image Milnor number equal to 1 
(this last a result already known when $n\geq p$). 
\end{abstract}

\section{Introduction}

In this paper we investigate the topology of the discriminant
of stable perturbations $f_t$ of multi-germs $f\!:\!(k^n,S)
\to (k^p,0)$  with $n \geq p-1$, where $S
\subset k^n$ is a finite set, and where $k=\RR$ or $\CC$. When $n=p-1$ 
`discriminant' of course means `image'.

When $k=\CC$, 
the discriminant $D(f_t)$ has the homotopy type of a wedge of $(p-1)$-spheres
(\cite{damon-mond},\cite{mondv}).
The number of these spheres is called the {\bf discriminant
Milnor number} when $n\geq p$ and the {\bf image Milnor number} when $n=p-1$,
and denoted $\mu_\Delta$ and $\mu_I$ respectively.
When $n \geq p$ and 
$(n,p)$ are in Mather's range of nice dimensions (\cite{mather6}),
it is known (\cite{damon-mond}) that $\mu_\Delta(f)$ and the ${\cal A}_e$-
codimension of $f$ satisfy the Milnor-Tjurina relation:
$$\mu_\Delta(f) \geq {\cal {A}}_e\mbox{-codimension} (f)$$
with equality if $f$ is weighted homogeneous in some coordinate system.
In case $n=p-1$, the same relation,
with $\mu_I$ in place of $\mu_\Delta$, is only known to hold when
$n=1$ (\cite{curves}) and $n=2$ (\cite{js},\cite{mondv}).
Nevertheless there is evidence that it holds 
in higher dimensions (see e.g. 
\cite{houston-kirk}):\\

\ni{\bf Conjecture I} This relation holds in all nice dimensions $(n,n+1)$.\\

Here we are concerned with this conjecture, and also with another:
suppose that
$g:(\RR^n,S)\to (\RR^p,0)$ is a real analytic map germ of finite
$\cal {A}$-codimension, with a stable perturbation $g_t$. Suppose also
that the complexification $g_{\CC,t}$
of $g_t$ 
is a stable perturbation of the complexification $g_\CC$ of $g$. We
say that $g_t$ is a {\bf good real perturbation} of $g$ 
if $\mbox{rank}\,H_{p-1}
(D(g_t);\Bbb{Z})=\mbox{rank}\,H_{p-1}(D(g_{\Bbb{C},t});\Bbb{Z})$ (in which
case the inclusion of real in complex induces an isomorphism on the vanishing
homology of the discriminant).\\

\ni{\bf Conjecture II} 
For every ${\cal A}_e$-codimension 1 equivalence class of 
map-germs in 
the nice dimensions, there exists a real form with a good real perturbation.
That is, the vanishing topology of all codimension 1 complex singularities
is `visible over $\Bbb{R}$'.\\

\ni For maps $\CC^2\to\CC^3$ there are five codimension 1 equivalence- classes 
(see Figure 1 on page 10);
for maps $\CC^3\to\CC^4$ there are eight, and for maps $\CC^4\to\CC^5$ there 
are eleven.\\
  
\ni Conjecture II is known to hold for mono-germs of maps 
$\CC^n\to\CC^p$ (with $n\geq p$ and
$(n,p)$ nice dimensions) of corank 1 (\cite{mond}). It also holds for 
(mono- and multi-)
germs of maps $\CC^2\to\CC^3$ (\cite{goryunov}; Goryunov's diagrams of 
good real perturbations are reproduced in Figure 1 below).  
Every real germ $\CC\to\CC^2$ has a 
good real perturbation (\cite{acampo},\cite{gusein}), but once $n>1$, 
map-germs $\CC^n\to\CC^{n+1}$ with good
real perturbations become the exception (\cite{marar-mond}).   

Our main results here provide  evidence for both conjectures.
We show\\
 
\ni{\bf Theorem \ref{inumber}} Every multi-germ $f:(\CC^n,S)\to
(\CC^{n+1},0)$ of corank
1 and ${\cal{A}}_e$-codimension 1 has $\mu_I(f)=1$.\\

\ni{\bf Theorem \ref{grp}} Every  ${\cal A}$-equivalence class of 
multi-germ $f:(\CC^n,S)\to(\CC^p,0)$ ($n\geq p-1, (n,p)$ nice dimensions) 
of corank
1 and ${\cal{A}}_e$-codimension 1 has a real form with a good real 
perturbation.\\

We prove both of these theorems first
for 'mono-germs' ($|S|=1$) (in Section~\ref{mono-germs}) and then by an 
inductive procedure which
constructs codimension 1 multi-germs from simpler ingredients. This procedure
yields an inductive classification of multi-germs of codimension 1. 
In Section \ref{multi-germs} we show that all codimension 1 
multi-germs can be constructed from codimension $1$ multi-germs
with fewer branches and in a lower dimension, and from trivial unfoldings of
Morse singularities (in case $n\geq p$) or immersions (in case $p=n+1$)
by means of three standard operations. These are {\it augmentation}, described
in Section \ref{augmentations}, 
and two {\it concatenation operations}, described in Section
\ref{concatenation}.\\ 

We feel that these operations, of augmentation and concatenation, are 
themselves of independent interest. They can be seen at work, generating
the lists of ${\cal A}_e$-codimension 1 germs from surfaces to 3 space, 
and from surfaces to surfaces, in Figures 1 and 2, on page 10. See also
Figure 3, on page 13.\\

We end this introduction with an elementary lemma which nevertheless
highlights an important property of codimension 1 germs.
\begin{lema}\label{stabver}
If $f:(\CC^n,S)\to(\CC^p,0)$ 
is a germ of ${\cal A}_e$-codimension 1, then any stable unfolding
of $f$ is ${\cal A}_e$-versal.
\end{lema}
\Proof Let $F(x,u)=(f_u(x),u)$ be a $d$-parameter unfolding of $f$.
By \cite{martinet} XV 2.1, $F$ is stable iff
$$T{\cal A}_ef+\OO_p\{\partial f_u/\partial u_1|_{u=0},\ldots,\partial 
f_u/\partial u_d|_{u=0}\}=\theta(f).$$
Since $T{\cal A}_ef$ is an $\OO_p$-module, we therefore cannot have 
$\partial f_u/\partial u_i|_{u=0}\in T{\cal A}_ef$ for all $i$.
Hence for some $i$, $T{\cal A}_ef+\CC\{\partial f_u/\partial u_i|_{u=0}\}
=\theta(f)$, and $F$ is versal.\eop

The results in this paper concerning maps $\CC^n\to\CC^{n+1}$, and the results
of Sections 2 and 5, 
were first proved in the Ph.D. thesis (\cite
{cooper}) of the
first author. 
\section{Augmentations}\label{augmentations}

Let $f\!:\!(\Bbb{C}^n,S) \to (\Bbb{C}^p,0)$\, be a multi-germ of ${\cal
{A}}_e$-codimension 1 where $S$ is a finite subset of $\Bbb{C}^n$. Let
$$\begin{array}{c}
F\!:\!(\Bbb {C} \times \Bbb {C}^n, \{0\} \times S) \to (\Bbb {C} \times \Bbb
{C}^p, (0,0))\\
(\lambda,x) \mapsto (\lambda,f_\lambda (x))
\end{array}$$
be an ${\cal {A}}_e$-versal unfolding of $f$. Define $A_F:(\Bbb{C}
\times \Bbb {C}^n, \{0\} \times S) \to (\Bbb {C} \times \Bbb {C}^p,
(0,0))$\, by
$A_F(\lambda,x)=(\lambda, f_{\lambda^2}(x))$.

\begin{prop}\label{indep1}

The $\cal {A}$-equivalence class of $A_F$  
is independent of the choice of miniversal unfolding $F$ of $f$. It depends
only on the $\cal A$-equivalence class of $f$.
\end{prop}
\Proof 
Let $F(t,x)=(t,f_t(x))$ and $G(s,x)=(s,g_s(x))$ be two 1-parameter versal
unfoldings of $f$. 
From the definition of versality it follows immediately
that there exist diffeomorphisms $\Phi(t,x)=(t,\phi_t(x))$ and $\Psi(t,y)
=(t,\psi_t(y))$ and a base-change diffeomorphism $\alpha:(\CC,0)\to(\CC,0)$
such that
$\alpha^\ast(F)(t,x)=\Psi\circ G\circ \Phi$ (where $\alpha^\ast(F)$ is the
unfolding $(t,x)\mapsto (t,f_{\alpha(t)}(x)$).
An easy calculation shows that there exists $\beta:(\CC,0)\to(\CC,0)$ (also
invertible) such that $\alpha(t^2)=\beta(t)^2$; now writing
$A_\Phi(t,x)=(t,\phi_{t^2}(x))$ and $A_\Psi(t,y)=(t,\psi_{t^2}(y))$ we have
$\beta^\ast(A_F)=A_\Psi\circ A_G\circ A_\Phi$. 

Equivalence of germs entails equivalence of their miniversal unfoldings,
so the second statement follows.
\eop

We shall write $Af$ for the $\cal {A}$-equivalence class of $A_F$.
We call $Af$ the {\bf augmentation} of $f$ and say that a multi-germ is an
augmentation if and only if it is the augmentation of some multi-germ $f$. A
multi-germ that is not an augmentation is called {\bf primitive}.

\begin{example}\label{goodpics}
{\em The five ${\cal {A}}_e$-codimension 1
multi-germs from $\Bbb {C}^2$ to $\Bbb {C}^3$
are:

I. $S_1$ (the birth of two umbrellas).

II. The non-transverse contact of two immersed sheets.

III. The intersection of three immersed sheets which are pairwise
transverse, but with
each one having first order tangency to the intersection of the other two.

IV. A cross-cap meeting an immersed plane.

V. A quadruple intersection.

IV and V are primitive.
I is the augmentation of the
cusp $t\mapsto (t^2,t^3)$, II is the augmentation of a tacnode
(two curves simply tangent at a point), which itself is the
augmentation of the map from two copies of $\Bbb {C}^0$ to $\Bbb {C}$
sending both
points to $0 \in \Bbb {C}$, and III is the augmentation of
three
lines
meeting pairwise transversely at a point. \\

\ni Pictures of the images of good real perturbations of these germs, showing
the process of augmentation, are shown in 
Figure 1 on page 10.
}
\end{example}

\begin{example}\label{betterpics} {\em The five ${\cal {A}}_e$-codimension 1
multi-germs from $\Bbb {C}^2$ to $\Bbb {C}^2$
are:

I. The lips: $(x,y)\mapsto (x,y^3+x^2y)$;

II. The swallowtail: $(x,y)\mapsto (x,y^4+xy)$;

III. The fold tacnode - a bi-germ consisting of two folds whose discriminant
curves have a simple tangency;

IV. The fold triple-point: a tri-germ consisting of three folds whose 
discriminants meet pairwise transversely at a point;

V. A bi-germ consisting of a fold and a cusp, with the discriminant of the 
fold transverse to the limiting tangent line to the discriminant of the cusp.\\

\ni I is the augmentation of $y\mapsto y^3$; II is primitive; III is the 
augmentation of the bi-germ consisting of the two branches $x\mapsto x^2$ and
$y\mapsto y^2$; IV and V are both primitive.

\ni Pictures of the discriminants of good real perturbations of these germs
are shown in Figure 2 on page 10.
}
\end{example}

\begin{teo}\label{cod}

$Af$ has ${\cal {A}}_e$-codimension 1.

\end{teo}
\Proof We use Damon's theory of ${\cal {K}}_V$-equivalence (see for
example~\cite{damon}).

The diagram
$$\begin{array}{ccc}
\Bbb {C}^{n+1} & \stackrel{F}{\rightarrow} & \Bbb {C}^{p+1}\\
\uparrow \mbox{id}& & \,\,\,\uparrow \gamma \\
\Bbb {C}^{n+1} & \stackrel{A_Ff}{\rightarrow} & \Bbb {C}^{p+1} \end{array}$$
where $\gamma(\delta,y)=(\delta^2,y)$, is a transverse fibre
square.
Therefore the ${\cal {A}}_e$-codimension of $A_Ff$ is equal to the ${\cal
{K}}_{D(F),e}$-codimension of $\gamma$ where $D(F)$ is the discriminant of
$F$.

But the diagram
$$\begin{array}{ccc}
\Bbb {C}^{n+1} & \stackrel{F}{\rightarrow} & \Bbb {C}^{p+1}\\
\,\,\,\uparrow i_1& & \,\,\,\uparrow i_2\\
\Bbb {C}^n & \stackrel{f}{\rightarrow} & \Bbb {C}^p \end{array}$$
where $i_1$ and $i_2$ are inclusions, 
is also a transverse fibre square.
So the
${\cal {K}}_{D(F),e}$-codimension of $i_2$ is equal to the ${\cal
{A}}_e$-codimension of $f$ and therefore is 1.

Because $i_2$ is a standard coordinate immersion, an easy calculation shows
$$N{\cal K}_{D(F),e}i_2=\frac{\OO_p}{d\lambda (i_2^\ast(\Der
(\log D(F))))}$$
(where $d\lambda (i_2^\ast(\Der(\log (D(F))))$ is the module consisting of the
coefficients of $\partial/\partial \lambda$ of the elements of $i_2^\ast(\Der
(\log D(F))$). A similar calculation gives
$$N{\cal K}_{D(F),e}\gamma=\frac{\OO_{p+1}}{d\lambda (\gamma^\ast(\Der
(\log D(F))))+(\delta)}$$
where the $(\delta)$ in the denominator comes from $\partial \gamma/\partial
\delta$. Clearly $N{\cal K}_{D(F),e}i_2$ and $N{\cal K}_{D(F),e}\gamma$
are isomorphic. \eop

\vspace{.2in}

The following result is a partial converse.

\begin{prop}\label{aug}

Suppose that $G(\lambda,x)=(\lambda,g_\lambda(x))$ is a one-parameter stable
unfolding of a multi-germ $g=g_0$ and suppose that 
$h(\lambda,x)=(\lambda,g_{\lambda^2}(x))$
has ${\cal {A}}_e$-codimension 1. Then $g$ has ${\cal
{A}}_e$-codimension 1 and $G$ is a versal unfolding of $g$. Thus $h$ is the
augmentation of $g$.

\end{prop}
\Proof It is immediate from the calculation in the proof of \ref{cod} that $g$
has ${\cal A}_e$-codimension 1. Versality of $G$ now follows by \ref{stabver}.
\eop

Given a stable map $f\!:\!(\Bbb {C}^n,S) \to (\Bbb {C}^p,0)$ let $Pf$
(the 'prism' on $f$) be the trivial 1-parameter unfolding of $f$. We
shall say that a map-germ is a {\bf prism} if it is $\cal A$-equivalent to
$Pg$ for some germ $g$.

An easy calculation with tangent spaces shows
\begin{prop}

Let $F(\lambda,x)=(\lambda,f_\lambda(x))$ be an
${\cal {A}}_e$-versal unfolding of
an ${\cal {A}}_e$-codimension 1 multi-germ $f$. Then
$G(\mu,\lambda,x)=(\mu,\lambda,f_{\lambda^2+\mu}(x))$
is an ${\cal {A}}_e$-versal
unfolding of $g=A_Ff$. \eop
\end{prop}

Since $G(\mu,\lambda,x)=(\mu,\lambda,f_{\lambda^2+\mu}(x))$ is an
unfolding of
$F(\mu,x)=(\mu,f_\mu(x))$ and $F$ is stable then $G$ is
$\cal {A}$-equivalent to $PF$.
Therefore if a multi-germ is an augmentation, its miniversal
unfolding is a prism. The converse is also true:

\begin{teo}\label{augment}

Let $g$ be a multi-germ of ${\cal {A}}_e$-codimension 1 and suppose that the
miniversal unfolding $G$ of $g$ is a prism. Then $g$ is an augmentation.

\end{teo}
\Proof There is a unique natural number
$\ell$ and a stable multi-germ $h$, unique up to $\cal A$-equivalence,
such that $G(\lambda,x)=(\lambda, g_\lambda(x))$ 
is $\cal {A}$-equivalent to $P^\ell h$ and
$h$ is not a prism.

We have the following commutative diagram
$$\begin{array}{ccc}
\Bbb {C}^n,S & \stackrel{g}{\rightarrow} & \Bbb {C}^p,0\\ \,\, \downarrow & &
\,\,\,\,\,\,\, \downarrow i\\ \Bbb {C} \times \Bbb {C}^n,\{0\} \times S &
\stackrel{(\lambda, g_\lambda(x))}{\rightarrow} & \Bbb {C} \times \Bbb
{C}^p,(0,0)\\ \,\,\,\,\,\,\, \downarrow \phi & & \,\,\,\,\,\,\, \downarrow
\psi\\
\Bbb {C}^\ell\times\Bbb {C}^{n+1-\ell},\{0\}\times S' & 
\stackrel{id_{\CC^\ell}\times h}{\rightarrow} & \Bbb {C}^\ell \times 
\Bbb {C}^{p+1-\ell},(0,0)\\ \,\,
\downarrow & &
\,\,\,\,\,\,\, \downarrow \pi\\ \Bbb {C}^{n+1-\ell},S' &
\stackrel{h}{\rightarrow} &
\Bbb {C}^{p+1-\ell},0 \end{array}$$
where $i$ is the standard inclusion, $\phi$ and $\psi$ are
diffeomorphisms,
$\pi$ is the natural projection and $S'$ is a subset of $\Bbb {C}^{n+1-\ell}$
of the
same cardinality
as $S$. Each of the three squares of the diagram is a transverse
fibre square so the outside rectangle is a transverse fibre square as
well. The
${\cal {A}}_e$-codimension of $g$ is equal to the ${\cal
{K}}_{D(h),e}$-codimension of $\pi\circ\psi\circ i$ where $D(h)$ is the
discriminant of $h$. Since $h$ is stable it is Thom transversal so any vector
field in $\Der(\log D(h))$ lifts, by 6.14 of~\cite{looijenga}. Since $h$
is not a
prism, $\Der(\log D(h))\subseteq m_{p+1-\ell} \theta(p+1-\ell)$. So,
$$T {\cal {K}}_{D(h),e}(\pi\circ\psi\circ i)\subseteq T {\cal {K}}_e
(\pi\circ\psi\circ i)$$
and the ${\cal {K}}_e$-codimension of $\pi\circ\psi\circ i$
is 0 or
1. It cannot be 0, as this would make $\pi\circ\psi\circ i$ a submersion and
$g$ stable. Therefore $\pi\circ\psi\circ i$ is a quadratic singularity,
$\cal {A}$-equivalent to
$$(y_1,\ldots,y_p) \stackrel{\gamma}{\mapsto} (y_1,\ldots,y_{p-\ell},
\sum_{i=p+1-\ell}^p y_i^2)$$

Let $\Phi$ and $\Psi$ be germs of diffeomorphisms such that
$\Psi\circ(\pi\circ\psi\circ i)=\gamma\circ\Phi$.

Let $\pi_{p+1-\ell}\!:\!\Bbb{C}^{p+1-\ell} \to \Bbb{C}$ be projection onto the last
coordinate. Then $d(\pi_{p+1-\ell}\circ\Psi\circ
(\pi\circ\psi\circ i)) (0)=0$ and
since $h$ is transverse to $\pi\circ\psi\circ i$,
$d(\pi_{p+1-\ell}\circ\Psi \circ
h)(S') \neq 0$. It follows that for $\lambda$ near 0,
$(\pi_{p+1-\ell}\circ\Psi \circ
h)^{-1}(\lambda) \cong \Bbb{C}^{n-\ell}$ and
$(\pi_{p+1-\ell}\circ\Psi)^{-1}(\lambda)
\cong \Bbb{C}^{p-\ell}$.

Define $h_\lambda=h|_{(\pi_{p+1-\ell}\circ\Psi\circ h)^{-1}(\lambda)} :
\Bbb{C}^{n-\ell}
\to \Bbb{C}^{p-\ell}$. Then $h(\lambda,x)=(\l, h_\lambda(x))$ is an unfolding
of $h_0$.
Since the outside rectangle of the above diagram is a transverse fibre
square, $g$
is $\cal {A}$-equivalent to the germ $(\l_1,\ldots,\l_\ell,x)\mapsto (\l_1,\ldots,\l_\ell,h_{\sum_{i=1}^\ell\lambda_i^2}(x))$.
Therefore, $g$ is an augmentation by Proposition~\ref{aug}. \eop

\section{Concatenation}\label{concatenation}
In this section we describe two basic operations, by which we 
``concatenate'' stable unfoldings of (multi-) germs to create new
multi-germs. There is no reason to require purity of dimension in multi-germs,
and we allow different branches to have domains of different dimension.
We therefore will not distinguish in our notation between image Milnor
number and discriminant Milnor number: both will be denoted $\mu_\Delta$.
In what follows it will be useful to use the notation $\{f,g\}$ for
the germ obtained by putting together germs $f$ and $g$ with the same target.\\

\ni Throughout this section we assume that we are in the nice dimensions;
thus, every stable unfolding $(f_\lambda(x),\lambda)$ 
of a germ $f_0$ is a ``stabilisation'',
in the sense that for almost all $\lambda$, $f_\lambda$ is stable.

\ni The first concatenation operation
is monic: from a multi-germ with $m$ branches we get a 
multi-germ with 
$m+1$ branches, in which the extra branch is a fold or an immersion.
\begin{teo}\label{conc1}
Let $f:(\CC^n,S)\to(\CC^p,0)$ be a map-germ of finite ${\cal A}_e$-codimension
with a stable unfolding $F$ on the single parameter $t$, 
let $0\leq k\in\ZZ$ and let $g:(\CC^p\times\CC^k,0)\to
(\CC^p\times\CC,0)$ be the fold map $(y,v)\mapsto (y,\sum_{j=1}^kv_j^2)$. 
Then 
\be
\item
$${\cal A}_e-{\mbox codim}(g^\ast(F))={\cal A}_e\mbox{-codim}
(f)= {\cal A}_e\mbox{-codim}(\{F,g\})$$
\item
$$\mu_{\Delta}(g^\ast(F))=\mu_\Delta(\{F,g\})=\mu_\Delta(f)$$
\item
both $g^\ast(F)$ and $\{F,g\}$ have 1-parameter stable unfoldings. 
\ee
\end{teo}
\Proof (1) Let $i:(\CC^p,0)\to(\CC^p\times\CC,0)$ be the standard inclusion
inducing $f$ from $F$. 
By Damon's theorem the ${\cal A}_e$-codimension of $f$ is
equal to the vector-space dimension of $N{\cal K}_{D(F),e}i:=
\theta(i)/ti(\theta_{\CC^p})+
i^\ast(\Der(\log D(F)))$. As $i$ is an immersion, 
projecting to the last component gives an
isomorphism $N{\cal K}_{D(F),e}i\simeq 
\OO_{\CC^p,0}/dt(i^\ast(\Der(\log D(F))))$. 
Again by Damon's theorem, the ${\cal A}_e$-codimension of $g^\ast(F)$ is equal
to the dimension of $N{\cal K}_{D(F),e}g$; 
since $tg(\theta_{\CC^p\times\CC^k})=
\sum_{\ell = 1}^p\OO_{\CC^p\times\CC^k}\cdot \partial/\partial y_\ell+
\sum_{j=1}^k \OO_{\CC^p\times\CC^k}\cdot v_j\partial/\partial t$, 
it follows, again by projecting to the last component, that
$$N{\cal K}_{D(F),e}g\simeq \OO_{\CC^p\times\CC^k,0}/(v_1,\cdots,v_k)
+dt(g^\ast(\Der(\log D(F))));$$ 
this in turn is isomorphic to $\OO_{\CC^p,0}/
dt(i^\ast(\Der(\log D(F))))$, and thus to $N{\cal K}_{D(F),e}i$.
This proves the first equality in (1). 

To prove the second equality in (1), we use the exact sequence
$$0\to \frac{\theta(g)}{tg(\theta_{\CC^p\times\CC^k})+\omega 
g(\Der(\log D(F)))}
\to 
N{\cal A}_e\{F,g\}
\to
N{\cal A}_eF
\to
0
$$
which results from the fact that $\Der(\log D(F))$ is the kernel of 
$\overline{\omega} F\!:\!
\theta_{\CC^p\times\CC}\!\to\! \theta(F)/tF(\theta_{\CC^n\times\CC})$. 
Since $F$ is stable, $N{\cal A}_e
\{F,g\}$ is isomorphic to $\theta(g)/tg(\theta_{\CC^p\times\CC^k})+\omega g
(\Der (\log D(F)))$. This in turn is isomorphic to 
 $\OO_{\CC^p\times\CC^k,0}/(v_1,\cdots,v_k)+dt(\omega g(\Der(\log D(F)))$,
by projection to the last component, and thus, evidently, to
$\OO_{\CC^p,0}/dt(i^\ast(\Der(\log D(F)))$, i.e. to $N{\cal K}_{D(F),e}i$.\\

(2) For $\lambda\neq 0$, the map
$g_\lambda$ defined by $g_\lambda(y,v)=(y,\sum v_j^2+\lambda)$ is 
logarithmically transverse to $D(F)$. 
Thus $g_\lambda^\ast(F)$ is a stable perturbation of $g^\ast(F)$. Its
discriminant is $g_\lambda^{-1}(D(F))$. There are now two
cases, $k>0$ and $k=0$. 

If $k>0$, $g_\lambda^{-1}(D(F))$ fibres over
$D(F)$ with typical fibre diffeomorphic to the Milnor fibre $X_g$ of $g$, and
contractible fibres over the points of $D(F)\cap D(g_\lambda)$. Since $D(F)$
itself is contractible, it follows that $g_\lambda^{-1}(D(F))$ is homotopy-
equivalent
to the space obtained from $D(F)\times X_g$ by gluing in a $k$-ball to each 
fibre over $D(F)\cap D(g_\lambda)$ to kill its homotopy. A Mayer-Vietoris
argument now shows that the rank of $H_{p+k-1}(g_\lambda^{-1}(D(F)))$ is equal
to the rank of $H_{p-1}(D(F)\cap D(g_\lambda))$. 
Since $\{F,g_\lambda\}$
is a stable perturbation of $\{F,g\}$,  a second Mayer-Vietoris argument
shows that $H_p(D(F)\cup D(g_\lambda))\simeq H_{p-1}(D(F)\cap D(g_\lambda))$;
thus
$$\mu_\Delta(g^\ast(F))=\mbox{rank}\ H_{p+k-1}(g_\lambda^{-1}(D(F)))
=\mbox{rank}\ H_p(D(F)\cup D(g_\lambda))=\mu_\Delta\{F,g\}.$$
The second equality of (2) follows from the fact that
$D(g_\lambda)=i_\lambda(\CC^p)=D(i_\lambda^\ast(F))$, 
where $i_\lambda:\CC^p\to\CC^p\times\CC$
is defined by $y\mapsto (y,\lambda)$. For $i_\lambda$
is logarithmically transverse to $D(F)$, and thus $i_\lambda^{-1}(D(F))$
(for $\lambda\neq 0$)
is the discriminant of a stable perturbation $i_\lambda^\ast(F)$ of $f$.

If $k=0$, the situation is much simpler:  $g_\lambda^{-1}(D(F))$
is diffeomorphic to $D(g_\lambda)\cap D(F)$, and the assertion is proved by
a similar Mayer-Vietoris argument.\\

(3) The unfolding $G=(g_\lambda,\lambda)$ of $g$ induces from 
$F\times {\mbox id}_{\CC}$ a stable unfolding 
of $g^\ast(F)$, since it
is logarithmically transverse to $D(F)\times\CC$. The unfolding
$\{F\times{\mbox id}_{\CC},G\}$ of $\{F,g\}$ is stable, since the analytic
stratum $\CC^p\times\CC\cdot(1,1)$ of $G$ is tranverse to the analytic
stratum of $F\times{\mbox id}_{\CC}$.
\eop

In particular, if the 
germ $f$ satisfies Conjecture 1, then so does $\{F,g\}$. In fact, our proof of
\ref{conc1} shows that the same goes for the existence of good real
perturbations (Conjecture 2, in the case of map-germs of codimension 1).
\begin{teo}\label{goodreal1}
If $f$ has a good real perturbation then so does 
$\{F,g\}$, and vice versa.
\end{teo}
\Proof Replace $\CC$ by $\RR$ everywhere in the topological part of
the proof of \ref{conc1}. The Mayer-Vietoris argument shows that
$\mbox{rank}\,H_p(D_{\RR}(F)\cup D_{\RR}(g_\lambda))=
\mbox{rank}\, H_{p-1}(D_{\RR}(F)\cap D_{\RR}(g_\lambda)) 
=\mbox{rank}\,H_{p-1}(D(f_t))$, 
so that if either side has, for $t>0$ or for $t<0$,
rank equal
to the rank of the homology of the complexification, then so,
by \ref{conc1}, does the other. \eop

\begin{teo}\label{indep2} Suppose that
the germ $f$ of Theorem \ref{conc1} has ${\cal A}_e$-
codimension 1. 
Then up to ${\cal A}$-equivalence, the bi-germ $h=\{F,g\}$ obtained is
independent of the choice of stable unfolding $F$.
\end{teo}
\Proof 
Any stable 1-parameter unfolding of $f$ is 
also ${\cal A}_e$-versal. Thus, given two such, $F'$ and $F''$, by the
semi-uniqueness
of mini-versal unfoldings there are a diffeomorphism $\alpha:(\CC,0)\to(\CC,0)$
and unfoldings of the identity 
$\phi:(\CC^n\times\CC,S\times\{0\})\to(\CC^n\times\CC,S\times\{0\})$
and $\Psi:(\CC^p\times\CC,\{0\}\times\{0\})\to(\CC^p\times\CC,\{0\}\times\{0\})$
such that 
$$\Psi\circ F'\circ\Phi=\alpha^\ast(F''),$$ 
where $\alpha^\ast(F'')$ is the unfolding 
$(x,\lambda)\mapsto (f''(x,\alpha(\lambda)),\lambda)$. This equality can be
rewritten
$$(1\times\alpha)\circ\Psi\circ F'\circ\Phi\circ(1\times\alpha^{-1})= F'',$$
and therefore to 
conclude that $\{F',g\}$ and $\{F'',g\}$ are $\cal A$-equivalent, it remains
only to show that we can find a diffeomorphism $\theta$ such that
$$(1\times\alpha)\circ\Psi\circ g\circ\theta = g.$$
In fact we construct $\theta^{-1}$. Since
$$(1\times\alpha)\circ\Psi\circ g(y,v)
=(\psi(y,\sum v_j^2),\alpha(\sum v_j^2)),$$
we look for a diffeomorphism $\beta:(\CC^k,0)\to(\CC^k,0)$
such that 
$\alpha(\sum v_j^2)=\sum (\beta_j(v_1,\ldots,v_k))^2$. This seems easiest to
do by working directly with power series; for example, when $k=2$, and 
assuming for ease of notation that $\alpha'(0)=1$, we can take
$$\beta(v_1.v_2)=(v_1(1+\alpha_2(v_1^2+2v_2^2)+
\alpha_3(v_1^4+3v_1^2v_2^2+3v_2^4)+\cdots)^{1/2},
v_2(1+\alpha_2v_2^2+\alpha_3v_2^4+\cdots)^{1/2}),$$ 
where the $\alpha_i$ are the coefficients of the Taylor series of $\alpha$.
Now we find that
$$(\psi(y,\sum v_j^2),\alpha(\sum v_j^2))= g(\psi(y,\sum v_j^2),\beta(v));$$
the right hand side of this equality is the composite of 
$g$ with a diffeomorphism of its domain, and so we are done.
\eop
When $f$ has ${\cal A}_e$-codimension 1, 
the germ $g^\ast(f)$ obtained by applying the procedure of
theorem \ref{conc1} is the k-fold augmentation of $f$, $A^kf$. It will
be useful to have a notation for the multi-germ $\{F,g\}$: we will denote it
by $C_k(f)$. 
Both $A^kf$ and $C_k(f)$ are well-defined as 
$\cal A$-equivalence classes, by \ref{indep1} and \ref{indep2}.

\begin{example}
{\em 
Let $f=\{f_1,f_2,f_3,f_4\}$ be the stable multi-germ parametrising the
union of the four coordinate hyperplanes $\{x_i=0\}$ in $\CC^4$ (in
descending order of $i$), and let
$g(x,y,z)=(x,y,z,z+y+x^k)$. Then by successive de-concatenation,
the codimension and image Milnor number of the 5-germ 
$\{f,g\}$ are equal to those
of the 4-germ $g^\ast(f)$ and the 3-germ $(g^\ast(f_1))^\ast(\{g^\ast(f_2),
g^\ast(f_3),g^\ast(f_4)\})$. The latter is equivalent to
$$\cases{x\mapsto (x,-x^k))\cr x\mapsto (x,0)\cr x
\mapsto (0,x)\cr}.$$
This has ${\cal A}_e$-codimension and image Milnor number equal to $k$ --- 
an $r$-branch parametrised curve-germ in the plane has
$\mu_I=\delta-r+1$. It also has a god real perturbation, shown here when
$k=4$.
\vskip 10pt
\epsfxsize 5cm
\centerline{\epsfbox{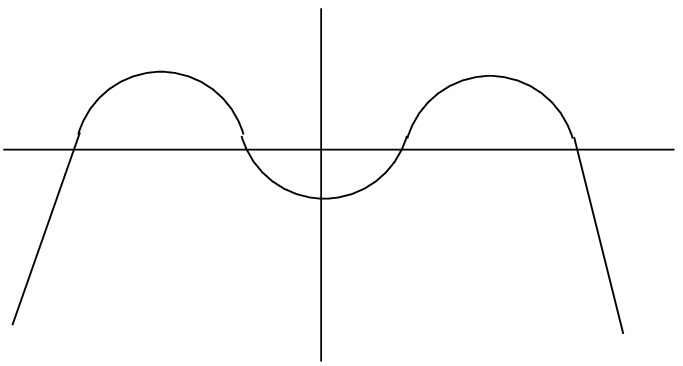}}
\vskip 10pt
}
\end{example}

%

%
%
\begin{example}\ 
{\em
\vskip 10pt
\epsfxsize 12cm
\centerline{\epsfbox{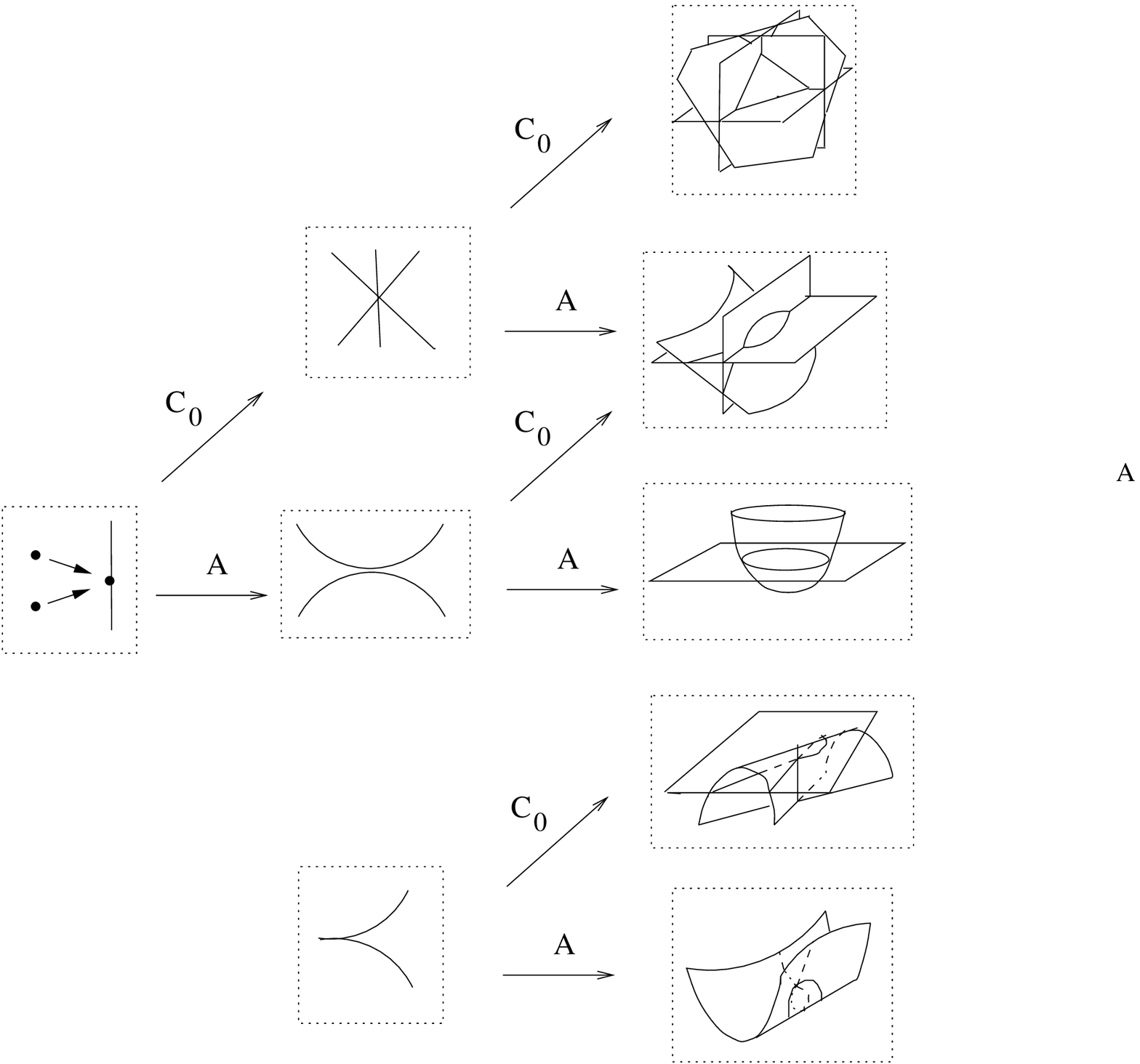}}
\begin{center}{{\em Figure 1: via $A=$Augmentation and $C_0=$Concatenation, 
the double-point and the cusp generate all the 
codimension 1 map-germs from 2-space to 3-space}}
\end{center}
}
\end{example}
\begin{example}\ 
{\em
\vskip 10pt
\epsfxsize 12cm
\centerline{\epsfbox{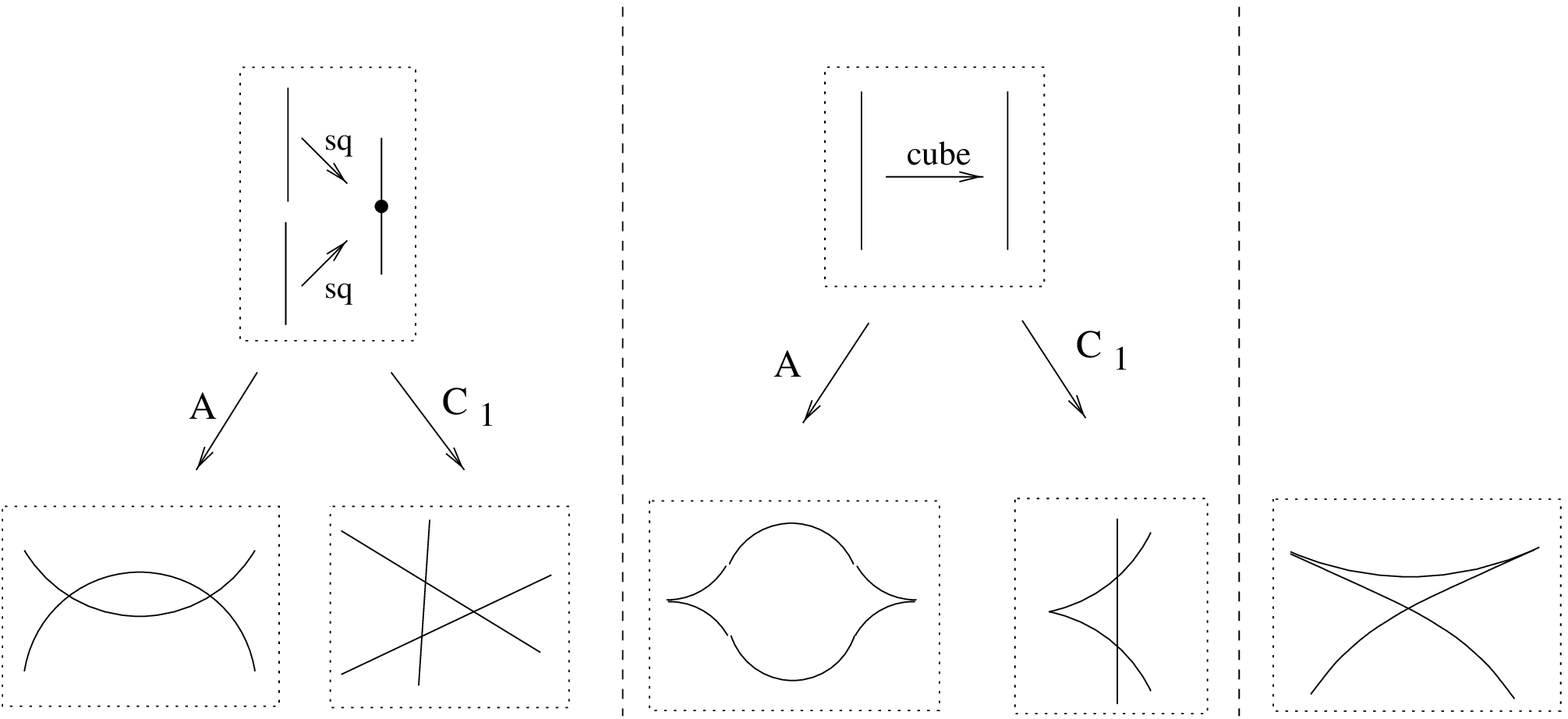}}
\begin{center}{{\em Figure 2: generation of 
codimension 1 germs of maps from the plane to
the plane}}
\end{center}
}\end{example}
\begin{example}{\em
The bi-germ consisting of a cross cap together with an
immersed plane transverse to the parametrisation of the cross-cap, and
making contact of degree $k$ with the double line in the cross-cap
(cf 7.5 in \cite{mondclass}, 3.3 in \cite{wik}) is obtained by applying $C_0$
to the germ $t\mapsto (t^2,t^{2k+1})$ parametrising the $k$-th order
cusp.}
\end{example}

The second type of concatenation is a binary operation: given 
germs $f_0\!:\!(\CC^m,S)\!\to\!(\CC^a,0)$ and $g_0:(\CC^n,T)\to(\CC^b,0)$ with 
1-parameter stable unfoldings $F$ and $G$, we form the multi-germ
$h$ essentially by putting together $\mbox{id}_{\CC^a}\times F$ and
$G\times\mbox{id}_{\CC^b}$ so that their analytic strata (see Section 5)
meet subtransversely
in $\CC^{a+b+1}$. 

\begin{teo}\label{conc2}
Suppose the two map-germs $F(y,s)=(f_s(y),s)$ and $G(x,s)=(g_s(x),s)$
are stable, and let 
$h$ be defined by
$$\cases{(X,y,u)\mapsto (X,f_u(y),u)\cr
(x,Y,u)\mapsto (g_u(x),Y,u)\cr}.
$$
Then provided ${\cal A}_e\mbox{-codim}\ (h) < \infty$,
we have
\be
\item
$${\cal A}_e\mbox{-codim}\ (h) \geq
{\cal A}_e\mbox{-codim}\ (f_0)\times {\cal A}_e\mbox{-codim}\ (g_0),$$
with equality if and only if
either $s\in ds(\Der(\log D(G)))$  or $t\in dt(\Der(\log D(F)))$;
\item
$h$ has a 1-parameter stable unfolding;
\item
$$\mu_\Delta(h)=\mu_\Delta(f_0)\times\mu_\Delta(g_0).$$\ee
\end{teo}
\Proof (1) and (2):
we compute the codimension of $h$ by Damon's
theorem. The multi-germ
$$H:\cases{(X,s,y,t)\mapsto (X,s,f_t(y),t)\cr
(x,s,Y,t)\mapsto (g_s(x),s,Y,t)\cr}
$$
is stable, as $\tau(F)\tr\tau(G)$, and after a change of coordinates can be 
seen as an unfolding of $h$
(which proves (2)).
Our map $h$ is induced from $H$ by 
$$i:\CC^a\times\CC^b\times\CC\to \CC^a\times\CC\times\CC^b\times\CC,$$
$$i(X,Y,u)=(X,u,Y,u).$$
The discriminant of $H$ is the ``product-union'' (Jim Damon's term) 
$$(D(G)\times\CC^b\times\CC)\bigcup (\CC^a\times\CC\times D(F)),$$
so if $\xi_0,\ldots,\xi_b$ generate $\Der(\log D(F))$ and
$\eta_0,\ldots,\eta_a$ generate $\Der(\log D(G))$ then (considering
the $\xi_i$ as belonging to $\theta(a+1+b+1/a+1)$ and the $\eta_i$ as
belonging to $\theta(a+1+b+1/b+1)$) we have
$$N{\cal K}_{D(H),e} i=\theta(i)/\langle
\partial/\partial X_i,
\partial/\partial Y_j,\partial/\partial s+\partial/\partial t\rangle
+\langle \xi_0,\ldots,\xi_b,\eta_0,\ldots,\eta_a\rangle.$$
Denote $ds(\Der(\log D(G)))$ and
$dt(\Der(\log D(F)))$ by $I$ and $J$ respectively.
By the map $(ds, dt)$, $N{\cal K}_{D(H),e}i$ projects isomorphically to
$M:=$
$$\frac{\OO_{a+b+1}\langle\partial/\partial s,\partial/\partial t\rangle}
{\langle \partial /\partial s+\partial/\partial t\rangle
+\langle \{\alpha(X,u)\partial/\partial s:\alpha(X,s)\in I\}\rangle
+\langle\{\beta(Y,u)\partial/\partial t:\beta(Y,t)\in J\}\rangle}.$$
As $f_0$ is induced from $F$ by $\gamma(y)=(y,0)$, and 
$g_0$ is induced from $G$ by $\sigma(x)=(x,0)$,
$$N{\cal A}_ef_0\simeq \theta(\gamma)/t\gamma(\theta_b)+
\gamma^\ast(\Der(\log D(F)))\stackrel{dt}{\simeq} \OO_b/\gamma^\ast(J)$$
and
$$N{\cal A}_eg_0\simeq \theta(\sigma)/t\sigma(\theta_a)+
\sigma^\ast(\Der(\log D(G)))\stackrel{ds}{\simeq} \OO_a/\sigma^\ast(I).$$
Now, suppose that $s\in I$. 
Then $M$ is isomorphic to $M_0:=$
$$
\frac{\OO_{a+b}\langle\partial/\partial s,\partial/\partial t\rangle}
{\langle \partial /\partial s+\partial/\partial t\rangle
+ \OO_{a+b}\sigma^\ast(I)\partial/\partial s\ +\ 
\OO_{a+b}\gamma^\ast(J)\partial/\partial t}$$
The reason that $M\simeq M_0$ is that
$u\partial/\partial s\in\{\alpha(X,u) 
\partial/\partial s):\alpha(X,t)\in
ds(\Der(\log D(G)))\}$ is in the denominator, and thus 
(since $\partial /\partial s+\partial/\partial t$ is in the denominator), 
so is $u\partial/\partial t$. \\ 

\ni Evidently, if  $t\in J$ then $M\simeq M_0$, by the
same argument. An easy argument shows that the converse is true: if
$M\simeq M_0$ then either $s\in I$ or $t\in J$.\\

\ni The module $M_0$ is itself isomorphic to 
$$\frac{\OO_{a+b}}{\sigma^\ast(I)+\gamma^\ast(J)}$$
via the map
$ds-dt$
$$ \alpha \partial/\partial s+\beta\partial/\partial t
\mapsto \alpha-\beta,$$
and finally, provided the left hand side is finite-dimensional,
$$\frac{\OO_{a+b}}{\sigma^\ast(I)+\gamma^\ast(J)}
\simeq \frac{\OO_a}{\sigma^\ast(I)}\otimes_{\CC}\frac{\OO_b}{\gamma^\ast(J)}.$$
This completes the proof of (1).

\ni (3) We postpone proof of this until Section \ref{topology}. 
\eop
\begin{obs}{\em
 Let $f_0:(\CC^n,S)\to(\CC^p,0)$ be a germ  
with a 1-parameter stable unfolding $F$, and suppose $n \geq p$
and $(n,p)$ are nice dimensions. 
Then the condition in the proposition,
that $t\in dt(\Der(\log D(F)))$, is equivalent to having $\mu_\Delta(f_0)=
{\cal A}_e-\mbox{codim}(f_0)$ - see \cite{GM}, Corollary 7.4. The proof uses
coherence of the Gauss-Manin connection.}
\end{obs}
Now suppose both $f_0$ and $g_0$ have ${\cal A}_e$-codimension 1.
By analogy with augmentation and the first type of concatenation, one would
expect the result of this second type of concatenation to be independent,
up to $\cal A$-equivalence, of
the choice of stable unfoldings $F$ and $G$. Somewhat surprisingly, this
is true over $\CC$ but false over $\RR$.
\begin{example}{\em
Let $f_0(y)=y^3,\ g_0(x)=x^3$, and take $F'(y,u)=(y^3+uy,u),\ 
F''(y,u)=(y^3-yu,u), G(x,u)=(x^3+ux,u)$. Then the multi-germs
$$h':\cases{(X,y,u)\mapsto (X,y^3+uy,u)\cr
(x,Y,u)\mapsto (x^3+ux,Y,u)\cr}$$
and
$$h'':\cases{(X,y,u)\mapsto (X,y^3-uy,u)\cr
(x,Y,u)\mapsto (x^3+ux,Y,u)\cr}$$
are not equivalent over $\RR$. The discriminant of $h'$ consists
of two components, each the product of a first-order cusp with a line,
and both ``opening downwards'' (in the direction of the negative $u$ axis).
This germ $h'$ does not have a good real perturbation. On the other hand,
in the germ $h''$ one cusp opens upwards and the other downwards, and 
$h''$ does have a good real perturbation, shown in Figure 3.  
\vskip 10pt
\epsfxsize 10cm
\centerline{\epsfbox{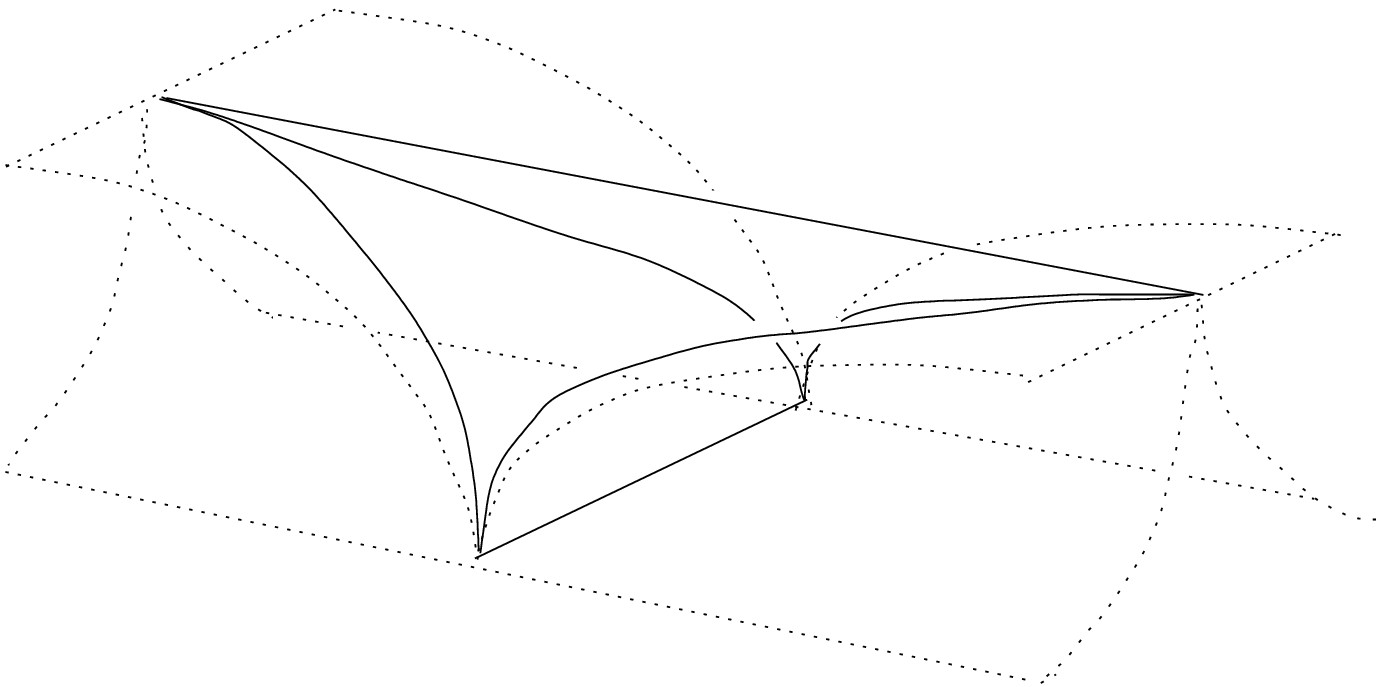}}
\begin{center}{\em Figure 3: Discriminant of a good real perturbation of
a binary concatenation of two
cubic functions}
\end{center}
\vskip 10pt

}
\end{example}
\begin{prop}\label{indep3}
Suppose that the germs $f_0$ and $g_0$ in Theorem \ref{conc2} both
have ${\cal A}_e$-codimen- sion 1. Then over $\CC$,
and up to $\cal A$-equivalence, the
germ $h$ is independent of choice of the 1-parameter
stable unfoldings $F$ and $G$.
\end{prop}
\Proof Suppose that $F'$ and $F''$ are 1-parameter stable unfoldings of 
$f_0$, and let $G$ be a 1-parameter stable unfolding of $g_0$. Applying
the concatenation operation, we obtain
multi-germs $h'$ and $h''$, the first using $F'$ and $G$, the second $F''$ 
and $G$. We wish to show that the two are $\cal A$-
equivalent. Let $h_\lambda$ be the linear interpolation between them:
$h_\lambda=(1-\lambda h')+\lambda h''$. We use a Mather-Yau type argument.\\

\ni{\bf Step 1} For no value of
$\lambda$ is the germ $h_\lambda$ stable.\\
For the analytic strata of its branches $\mbox{id}_{\CC^n}\times
F_\lambda$ and $G$ always meet at $0\in \CC^a\times\CC^b\times\CC$,
and always have dimensions whose sum is less than $a+b+1$, unless
for some value of $\lambda$ $F_\lambda$ is a trivial unfolding of 
$f_0$. In the latter case $F_\lambda$ itself is not stable, so that
once again $h_\lambda$ cannot be stable. 
It also follows that for those
$\lambda$ such that $h_\lambda$ has ${\cal A}_e$-codimension 1,
$T{\cal A}\,h_\lambda=T{\cal K}\,h_\lambda$. 
\\

\ni{\bf Step 2} The set of points $\{\lambda\in\CC:\ {\cal A}_e-\mbox
{codimension}(h_\lambda)>1\}$ is Zariski-closed in $\CC$, so
that its complement,  $\Lambda_1:=\{\lambda\in\CC:\ {\cal A}_e-\mbox
{codimension}(h_\lambda)=1\}$, is Zariski-open, and connected.
\\

\ni{\bf Step 3} Choose an integer $k$ such that in the 
appropriate multi-jet space $_rJ^k(X,Y)$, the $J^k{\cal A}$-orbit
of the $k$-jet of every codimension 1 germ coincides with the set of 
$k$-jets of its $\cal A$-orbit. We use Mather's Lemma (\cite{mather}, 3.1) 
to show that
the set $J^k\Lambda_1:=\{j^kh_\lambda:\lambda\in\Lambda_1\}$ lies in a 
single $J^k{\cal A}$-orbit, from which the proposition follows. 
It is necessary to check only that 
$T_\sigma J^k\Lambda_1\subset TJ^k{\cal A}\sigma$ for all 
$\sigma\in \Lambda_1$. But $J^k\Lambda_1$ lies in a single contact
orbit, and for each $\lambda\in \Lambda_1$, the ${\cal A}$-tangent space
of $h_\lambda$ is equal to its contact tangent space. It follows that
$T_\sigma J^k\Lambda_1\subset TJ^k{\cal A}\sigma$ for all
$\sigma\in \Lambda_1$, as required.
\eop
The argument of this proof in fact proves the following result, which we
will use later:
\begin{lema}\label{classifier}In any
given (complex)
contact class there is at most one open $\cal A$-orbit. \end{lema}
\eop
In the light of \ref{indep3}, we will refer to the ${\cal A}$-equivalence
class of multi-germ
obtained from codimension 1 multi-germs $f_0$ and $g_0$ by this binary 
concatenation operation 
as $B(f,_0,g_0)$.\\

\ni {\bf Question}\quad How many different $\cal A$-equivalence classes 
of germs $h$ over $\RR$
can the different choices of miniversal unfoldings $F,G$ of $f_0$ and $g_0$
give rise to? \\

\ni Our final result here is 
\begin{prop}\label{goodreal2}
If the germs $f_0$ and $g_0$ both have good real perturbations, then
so does $B(f_0,g_0)$.
\end{prop}
The proof will be given in Section \ref{topology}.
\begin{obs} {\em
It would be interesting to understand the effect on monodromy groups
of augmentation and concatenation. There is a ``natural''
choice of 1-parameter stable unfolding of $Af_0$,
$C_k(f_0)$ and of $B(f_0,g_0)$, reflecting the choice of stable unfolding
used in their construction. Presumably
the monodromy action in the case of $B(f_0,g_0)$ is the tensor product
of the monodromy action in the chosen 1-parameter unfoldings $F$ and $G$,
as in the classical Thom-Sebastiani theorem.
}
\end{obs}
\section{${\cal A}_e$-codimension 1 germs $(\CC^n,0)\to(\CC^{n+1},0)$}
\label{mono-germs}
In this section we first classify ${\cal A}_e$-codimension 1 mono-germs
and then show that each has image Milnor number 1. 
The argument runs roughly as follows: let
$$D^k(f)=\mbox{closure}
\{(x_1,\ldots,x_k)\in (\CC^n,S)^k|x_i\neq x_j \,\mbox{for}\
i\neq j, f(x_i)=f(x_j)\forall
i,j\};$$
then by results of \cite{mararmond},
$f$ is stable if and only if $D^k(f)$,
is smooth of dimension $n-k+1$, 
for $2\leq k\leq n+1$,and $f$ has finite ${\cal A}_e$-codimension 
if and only if each $D^k(f)$ is an isolated complete intersection
singularity of dimension $n-k+1$, again for $0\leq k\leq n+1$. Moreover,
if $f_t$ is a stable perturbation of $f$, then $D^k(f_t)$ is a Milnor
fibre of $D^k(f)$.
There is an obvious symmetric group action on $D^k(f)$, permuting the copies
of $(\CC^n,S)$, and in fact
a spectral sequence (\cite{goryunovmond}) computes the homology 
of the imageof $f_t$ from the $S_k$-alternating part of the 
homology of $D^k(f_t)$.It turns out that if $f$ has 
${\cal A}_e$-codimension 1, then just one of
the $D^k(f)$ is singular, and in fact has a Morse singularity. Since
the symmetric group action on the Jacobian algebra is therefore
trivial, from a theorem of Orlik-Solomon and Wall it follows that
the vanishing homology of $D^k(f_t)$ is alternating, and thus by
the spectral sequence the image Milnor number is 1.
The symmetry of $D^k(f)$ also accounts for the 
existence of a good real 
perturbation. Essentially, the point is that an $S_k$-invariant 
Morse function in 
$k$ real 
variables is either a sum of squares or the negative of a sum of squares.\\

Now we proceed with the classification. Let $\ell > 0$, take coordinates
$(u_1,\cdots,u_{\ell -1},v_1,\cdots,v_{\ell -1},x)$ on $\Bbb{C}^{2\ell
-1}$, and
define a map \\
$f^\ell:(\Bbb{C}^{2\ell -1},0)\to(\Bbb{C}^{2\ell},0)$ by
$$f^\ell(u,v,x)=(u,v,x^{\ell+1}+
\sum_{i=1}^{\ell -1}u_ix^i,x^{\ell+2}+\sum_{i=1}^{\ell -1}v_ix^i).$$
\begin{lema}\label{Dk}
The map-germ $f^\ell$ just described has ${\cal A}_e$-codimension 1, and the
following property:\\
$(\ast)$:
$D^k(f^\ell)$ is smooth for $2\leq k\leq {\ell}$, $D^{\ell+1}(f^\ell)$ has a Morse
singularity, and $D^k(f^\ell)$ is empty for $k>\ell+1$. 
\end{lema}
\Proof
Recall from \cite{mararmond}
2.1 the determinantal equations $h^k_{j,i}$ of $D^k(f^\ell)$: $$
h^k_{j,i}=\frac{\left|\begin{array}{cccccccc}
1&x_1&\cdots&x_1^{i-1}&f^\ell_j(u,v,x_1)&x_1^{i+1}&\cdots&x_1^{k-1}\\
&&&&\cdot&&&\\
&&&&\cdot&&&\\
&&&&\cdot&&&\\
1&x_k&\cdots&x_k^{i-1}&f^\ell_j(u,v,x_k)&x_1^{i+1}&\cdots&x_k^{k-1}
\end{array}\right|
}{\mbox{vdM}}
$$
for $1\leq i\leq k-1$ and $2\ell - 1\leq j\leq 2\ell$,
where vdM is the van der Monde determinant of
$x_1,\cdots,x_k$, and $f^\ell_j$ is the
$j$'th component of $f^\ell$. An easy calculation
shows
$$h^k_{2\ell -1,i}=u_i+O(2)\hspace{.5in}\mbox{for $i=2,\cdots,\ell-1$}$$
$$h^k_{2\ell,i}=v_i+O(2)\hspace{.5in}\mbox{for $i=2,\cdots,\ell-1$}$$ so that
$D^k(f^\ell)$ is smooth for $2\leq k\leq \ell$; moreover $$h^{\ell +1}_{2\ell
-1,\ell}=x_1+\cdots+x_{\ell+1}$$ and
$$h^{\ell +1}_{2\ell,\ell}= \sum_{i,j=1}^{\ell +1}x_ix_j.$$ We may take the
$h^{\ell +1}_{j,i}$ for $2\leq i\leq \ell-1$, together with
$x_1,\cdots,x_{\ell
+1}$, as coordinates; then $D^{\ell+1}(f)$ is embedded in
$x_1,\cdots,x_{\ell+1}$-space with equations $h^{\ell +1}_{2\ell -1,\ell}$
and
$h^{\ell +1}_{2\ell,\ell}$. Now $h^{\ell +1}_{2\ell -1,\ell}$ is
non-singular, and
$$h^{\ell +1}_{2\ell,\ell}-\frac{1}{2}(h^{\ell +1}_{2\ell -1,\ell})^2
=\sum_{i=1}^{\ell +1}x_i^2,$$
so $D^{\ell +1}(f^\ell)$ has a Morse singularity at the origin.\\

Calculation of the ${\cal A}_e$-codimension of $f^\ell$ is straightforward; it
may easily
be checked using nothing more than Nakayama's Lemma that $$T{\cal
A}_ef^\ell=\theta(f)\setminus \{x^\ell\p/\p Y_2, x^{\ell-1}\p/\p v_1, \cdots,
x\p/\p
v_{\ell-1}\}+$$
$$+\langle x^{\ell-1}\p/\p v_1+x^\ell\p/\p Y_2,\cdots, x\p/\p
v_{\ell-1}+x^\ell\p/\p Y_2\rangle.$$ The
calculation is carried out in detail in
\cite{cooper}. \eop
Note that since $f^\ell$ has ${\cal A}_e$-codimension 1, its $\cal A$-orbit
is open
in its
$\cal K$-orbit. Note also that from the expression for $T{\cal A}_ef^\ell$ given
in the proof, it follows that the stable germ
$$F(\lambda,u,v,x)=(\lambda,u,v,x^{\ell+1}+ \sum_{i=1}^{\ell
-1}u_ix^i,x^{\ell+2}+\sum_{i=1}^{\ell -1}v_ix^i+ \lambda x^\ell)$$
is an ${\cal A}_e$-versal
unfolding of $f^\ell$. 

Since for corank 1 germs $(\CC^n,0)\to(\CC^{n+1},0)$ the multiplicity
determines the contact class, it follows from Lemma \ref{classifier} that
we have
\begin{corl}\label{classified} If
$f:(\Bbb{C}^{2\ell -1},0)\to (\Bbb{C}^{2\ell},0)$
has corank 1,
multiplicity $\ell+1$ and ${\cal A}_e$-codimen-\\
sion 1, then $f$ is
$\cal A$-equivalent to the germ $f^\ell$ of Lemma ~\ref{Dk}. \eop
\end{corl}
\begin{prop}\label{classification}
If $f:(\Bbb{C}^n,0)\to(\Bbb{C}^{n+1},0)$ has corank $1$, multiplicity
$\ell+1$ and ${\cal A}_e$-codimension 1 then it is equivalent to
$$f^\ell_q:(u,v,w,x)\mapsto (u,v,w,x^{\ell+1}+\sum_{i=1}^{\ell-1}u_ix^i,
x^{\ell+2}+\sum_{i=1}^{\ell-1}v_ix^i+q(w)x^\ell)$$
where $q$ is a non-degenerate
quadratic form. \end{prop}
\Proof Note that $f^\ell_q$ is (over $\CC$) equivalent to the $k$-fold
augmentation $A^kf^\ell$, where $k=n-2\ell$.
The hypothesis forces $n\geq 2\ell-1$, since the minimal target
dimension
of a stable corank 1 germ of multiplicity $\ell+1$ is $2\ell+1$. Since $f$
has
${\cal A}_e$-codimension 1, its versal unfolding
$G:(\Bbb{C}^n\times\Bbb{C},0)\to(\Bbb{C}^{n+1}\times\Bbb{C},0)$ is an
$n-2\ell+1$-fold
prism
on a minimal stable map-germ of multiplicity $\ell+1$.
From this it follows by Theorem~\ref{augment} that
$f$ is equivalent to an
$n-2\ell+1$-fold augmentation of an ${\cal A}_e$-codimension 1 germ
$f_0:(\Bbb{C}^{2\ell-1},0)\to(\Bbb{C}^{2\ell},0)$ of multiplicity $\ell+1$
and corank
1. By the previous corollary, $f_0$ is equivalent to the germ $f^\ell$ of
~\ref{Dk};
since the germ $F$ described after ~\ref{Dk} is a versal unfolding of
$f^\ell$, $f$
is equivalent to the germ obtained by replacing the unfolding term $\lambda
x^\ell$ in the last component of $F$ by $q(w)x^\ell$, where $q$ is a
non-degenerate quadratic form in new variables $w_i$, as required.
\eop
\begin{prop}\label{rep}
If $f:(\Bbb{C}^n,0)\to(\Bbb{C}^{n+1},0)$ has corank 1 and ${\cal
A}_e$-codimension 1
then $\mu_I(f)=1$, and there is a real form with a good real perturbation.
\end{prop}
\Proof Let $f_t$ be a stable perturbation of $f$, with image $Y_t$. By
\cite{goryunovmond} Theorem 2.5,
\begin{equation}\label{homim}
H^n(Y_t;\Bbb{Q})\simeq \oplus_k \mbox{Alt}_kH^{n-k+1}(D^k(f_t);\Bbb{Q})
\end{equation}
where $\mbox{Alt}_kH^{n-k+1}(D^k(f_t);\Bbb{Q})$ means the subspace of
$H^{n-k+1}(D^k(f_t);\Bbb{Q})$ on which the symmetric group
$S_k$ acts by its sign
representation. Now $D^k(f_t)$ is a Milnor fibre of $D^k(f)$; since $f$ has
property $(\ast)$, (\ref{homim}) reduces to $$H^n(Y_t;\Bbb{Q})\simeq
\mbox{Alt}_{\ell +1}H^{n-\ell}(D^{\ell +1}(f_t);\Bbb{Q}).$$
As $D^{\ell +1}(f)$
has a Morse singularity,
$H^{n-\ell}(D^{\ell +1}(f_t);\Bbb{Q})\simeq \Bbb{Q}$; it
remains to show that the representation of
$S_{\ell +1}$ on $H^{n-\ell}(D^{\ell
+1}(f_t);\Bbb{Q})$ is the sign representation. This can easily be seen by an
explicit calculation with the normal form given; but there is another
argument which explains better why it is true. As
$D^{\ell +1}(f)$ is an $S_{\ell+1}$-invariant hypersurface
singularity, by the theorem
of Orlik-Solomon and Wall (\cite{orliksolomon},\cite{wall}),
$$H^{n-\ell}(D^{\ell
+1}(f_t);\Bbb{Q})\simeq \mbox{Jac}_{D^{\ell +1}(f)}
\otimes_{\Bbb{Q}}\wedge^{\ell}(V)^\ast$$ as $S_{\ell +1}$ representations,
where
$V$ is an $S_{\ell +1}$-invariant smooth space containing $D^{\ell +1}(f)$
as a
hypersurface, and $\mbox{Jac}_{D^{\ell +1}(f)}$ is the Jacobian algebra of
$D^{\ell +1}(f)$. Since $D^{\ell +1}(f)$ is Morse, its Jacobian algebra is a
trivial $1$-dimensional representation of $S_{\ell +1}$, so
$H^{n-\ell}(D^{\ell
+1}(f_t);\Bbb{Q})\simeq \wedge^{\ell}(V)^\ast$. In fact we take
$V=D^{\ell+1}(G)$
where $G$ is a 1-parameter stable unfolding of $f$;
as noted above, $G$ is right-left equivalent to a suspension of $F$, and in
particular the $S_{\ell+1}$-action on $D^{\ell+1}(G)$ is equivalent to a
trivial
extension of the standard Weyl action $A_{\ell}$, in which $S_{\ell+1}$
acts on
$\{(x_1,\cdots,x_{\ell+1}):\sum_ix_i=0\}$ by permuting coordinates. Hence
$\wedge^{\ell}(V)^\ast$ is just the sign representation of $S_{\ell+1}$
and (as
vector spaces)
$$H^{2\ell -1}(Y_t;\Bbb{Q})\simeq
\mbox{Alt}_{\ell+1}H^{n-\ell}(D^{\ell+1}(f_t);\Bbb{Q})=\Bbb{Q}$$ so that
$\mu_I(f)=1$.\\

In the real case, we apply (~\ref{homim}) to a real stable perturbation
$f_{t,\Bbb{R}}$ of $f$, replacing $D^k(f_t)$ by $D^k(f_{\Bbb{R},t})$.
Consider
first the case $n=2\ell-1$, so $f$ is equivalent to the germ $f^\ell$ of
~\ref{Dk}.
Let $f_{\Bbb{R},t}^C$ be a stable perturbation. Evidently
$D^k(f_{\Bbb{R},t}^C)$
is contractible for $2\leq k < \ell+1$, and $D^{\ell+1}(f_{\Bbb{R},t})$ is
a real
Milnor fibre of a ${\ell-1}$-dimensional Morse singularity; hence it is a
$p$-sphere for some $p$ between $-1$ and $\ell -1$. We have to show that
either
for $t >0$ or $t<0$ it is an $\ell-1$-sphere. This follows from the fact that
$D^{\ell+1}(f^\ell)$ has a Morse singularity and an $S_{\ell+1}$-invariant
defining
equation, in a space in which the representation of $S_{\ell+1}$ is
equivalent to
the Weyl representation $A_{\ell}$ described above. Since the
representation is
irreducible, the stable manifold and unstable manifold of the gradient
flow must
be equal to $0$ and $V$ or $V$ and $0$ respectively, and any
$S_{\ell+1}$-invariant quadratic form must have index $0$ or
$\ell$. Since the versal unfolding $F$ of $f^\ell$ is a stable map,
$D^{\ell+1}(F)$
is smooth, and thus projection to the parameter space cuts out distinct real
Milnor fibres for $t>0$ and $t<0$. Hence at least one of these is an
$\ell$-sphere. Inclusion $D^{\ell+1}(f_{\Bbb{R},t}^\ell)\hookrightarrow
D^{\ell+1}(f_t^\ell)$
then induces an $S_{\ell+1}$-equivariant homotopy equivalence, so that the
representation of $S_{\ell+1}$ on
$H^{\ell-1}(D^{\ell+1}(f^\ell_{\Bbb{R},t}))$ is
once again the sign representation.

In the general case, let $f_{q\ \Bbb{R},t}^\ell$ be a stable
perturbation of $f_q^\ell$.
By taking $q=\sum_i w_i^2$, then $D^{\ell+1}(f_t)$ is an
$\ell-1+d$-dimensional sphere,
where $d=n-2\ell+1$ is the number of $w$-variables in the expression for
$f^\ell_q$
in ~\ref{classification}. In fact $D^{\ell+1}(f_{q\ \Bbb{R},t}^\ell)$ is the
join of
$D^{\ell+1}(_{\Bbb{R},t})$ and $q^{-1}(t)$, and the
representation of $S_{\ell+1}$ on its cohomology is just the sign
representation
as before.
\eop
\begin{obs} {\em The argument just used shows that if
$f:(\Bbb{C}^n,0)\to(\Bbb{C}^{n+1},0)$ has corank 1 and multiplicity $\ell+1$,
and has
a 1-parameter stable unfolding $F$, and if $D^k(f)$ is singular, then
$\mu_I(f)\geq {\ell+2-k}$. For from the fact that $D^k(f)$ is singular it
follows
that $D^j(f)$ is singular, for $k\leq j\leq \ell+1$. As $D^j(f)$ is a
hypersurface
in the smooth space $D^j(F)$, the argument used above can be applied. The
Jacobian
algebra of each singular $D^k(f)$ has $S_k$-invariant subspace of
dimension at
least 1 (since the constants form a 1-dimensional trivial representation),
and
hence by the theorem of Wall (rather than the earlier result of
Orlik-Solomon,
which applies only to weighted homogeneneous hypersurface singularities) the
alternating part of the middle homology of the Milnor fibre $D^k(f_t)$ has
rank at
least 1. The conclusion then follows by (\ref{homim})}.
\end{obs}
\section{${\cal {A}}_e$-codimension 1 multi-germs}\label{multi-germs}
In this section we show 
that in the nice dimensions all ${\cal A}_e$-codimension 1 multi-germs 
can be constructed by concatenation
and augmentation, beginning with stable germs and with 
primitive ${\cal A}_e$-codimension 1 mono-germs.

Submersive branches of multi-germs 
play a trivial role in classification and deformation theory, and 
we will ignore them in what follows.
In particular ``a multi-germ with $k$ branches" means a multi-germ with $k$
non-submersive branches.

For a multi-germ
$f\!:\!(\Bbb {C}^n,S) \to (\Bbb {C}^p,0)$ with branches $f^{(1)},
\ldots, f^{(s)}$, define
$$\tau (f)=
ev_0[(\omega f)^{-1}\{f^* m_p \theta (f)+tf(\theta (n)_S)
\}]$$
where $ev_0:\theta (p) \to T_0 \Bbb {C}^p$ is evaluation
at $0$, and
$$\tau' (f)=ev_0[(\omega f)^{-1}\{tf(\theta (n)_S)\}]$$

In fact $\tau' (f)=ev_0(\Der(\log D(f)))$ where
$D(f)$ is the discriminant (or image) of $f$.

The following result is due to Mather~\cite{mather}.

\begin{prop}\label{regular}

The multi-germ $f$ is stable if and only if each $f^{(i)}$ is stable and $\tau
(f^{(1)}),\\
\ldots,\tau (f^{(s)})$
have regular intersection with respect to $T_0
\Bbb {C}^p$.
Moreover, in this case 
$\tau (f)=\\
\cap_i\tau (f^{(i)})$.
\eop
\end{prop}

We now investigate the geometrical significance of $\tau'$.

\begin{lema}\label{tau}

If $f\!:\!(\Bbb {C}^n,S) \to (\Bbb {C}^p,0)$ is stable then $\tau
(f)=\tau' (f)$.
\eop
\end{lema}

\begin{lema} If $f=\mbox{id}_{\CC^m}\times g$ (i.e. $f=P^mg$) then
$\tau'(f)=T_0\CC^m\oplus\tau'(g)$.
\eop
\end{lema}

\begin{lema}\label{tauprism}
If $\dim_{\Bbb{C}}\tau'(f)=m$, then there is a germ $g$, not a prism,
such that $f\sim_{\AA}P^mg$. Moreover, if $\phi$ and $\psi$ are diffeomorphisms such that
$f\circ\phi=\psi\circ(\mbox{id}_{\CC^m}\times g)$, 
then $\tau'(f)=d\psi_0(T_0\CC^m\times\{0\})$.
\end{lema}
\Proof Suppose $tf(\xi)=\omega f(\eta)$.
If $\eta(0)\neq 0$ then also $\xi(s)\neq 0$ for $s\in S$, and
the orbits of $\xi$ and $\eta$ can be incorporated as coordinate lines
into new
coordinate systems on $\Bbb{C}^n,S$ and $\Bbb{C}^p,0$;
now the lemma just reduces to the
Thom-Levine Lemma (see e.g.\cite{duplessis}),
and $f\sim_{\AA}Pg_1$ for some germ $g_1$. Now apply the same procedure
to $g_1$. After $m$ iterations, we arrive eventually at a $g$ with
$\tau'(g)=0$, which is therefore not a prism.
\eop

\begin{prop}\label{anastratum}

If $f\!:\!(\Bbb {C}^n,S) \to (\Bbb {C}^p,0)$ and
$g\!:\!(\Bbb {C}^m,T) \to (\Bbb {C}^q,0)$ are multi-germs neither of which
are prisms and if $P^kf$ is $\cal
{A}$-equivalent to $P^\ell g$ then $|S|=|T|$, $n=m$, $p=q$, $k=\ell$ and $f$ is
$\cal
{A}$-equivalent to $g$.
Furthermore, if the $\cal {A}$-equivalence between $P^kf$ and $P^\ell g$
is given by
diffeomorphisms $\phi$ and $\psi$  as in the following diagram then 
$\psi(\Bbb {C}^k \times
\{0\})=\Bbb {C}^\ell \times \{0\}$
$$\begin{array}{lcccl}
\hskip 1.2in
&\Bbb {C}^k \times \Bbb {C}^n,\{0\} \times S & \stackrel{id_{\Bbb {C}^k}
\times
f}{\rightarrow} & \Bbb {C}^k \times \Bbb {C}^p,(0,0)&\\
\hskip 1.2in&
\,\,\,\,\,\,\,\downarrow
\phi & & \,\,\,\,\,\,\,\downarrow \psi&\\ 
\hskip 1.2in&
\Bbb {C}^\ell \times \Bbb
{C}^m,\{0\} \times
T & \stackrel{id_{\Bbb {C}^\ell} \times g}
{\rightarrow} & \Bbb {C}^\ell \times \Bbb
{C}^q,(0,0)&\hskip 1.7in \Box \end{array}
$$
\end{prop}

Given a multi-germ $f$, by Proposition~\ref{anastratum} there is
a well defined maximal sub-manifold of the target along which $f$ is trivial
(i.e. a prism). It is known as
the {\bf analytic stratum} of $f$, and coincides with the set-germ  of points
$y\in\CC^p,0$ such that the germ $f:(\CC^n,f^{-1}(y)\cap C_f)\to(\CC^p,y)$ is
$\AA$-equivalent to $f:(\CC^n,S)\to(\CC^p,0)$.
Moreover, $\tau'(f)$ is the
tangent space at 0 to the analytic
stratum of $f$.

\begin{prop}\label{tauspace}

Let $f\!:\!(\Bbb {C}^n,S) \to (\Bbb {C}^p,0)$ and
$g\!:\!(\Bbb {C}^n,T) \to (\Bbb
{C}^p,0)$ be multi-germs, and suppose that
$h=\{f,g\}$ has $\AA_e$-codimension 1.
Let $\psi$ be a germ of a 1-parameter family of maps 
$(\Bbb {C}^p,0)\to \CC^p$ such that $\psi_0= id_{\Bbb {C}^p}$ and
$$ev_0(\frac{d\psi_t}{dt}|_{t=0}) \not\in \tau'(f)+\tau'(g),$$
and write $G(\lambda,x)=(\lambda,(\psi_\lambda \circ g)(x))$.
Then $H:\ =\{\mbox{id}_{\CC}\times f,G\}$
is a versal unfolding of $h$.
\end{prop}
\Proof Write $H(\lambda,x)=(\lambda, h_\lambda(x))$.
If
$v=\frac{dh_\lambda}{d\lambda}|_{\lambda=0} \in T{\cal {A}}_e h$, then
$v=th(\xi)+\omega h(\eta)$ for some $\xi \in \theta (n)_{S\cup T}$ and
$\eta \in \theta (p)$.
It follows that $\omega f(\eta)= tf(-\xi)$ and $tg(-\xi)=\omega
g(\eta - \frac{d\psi_t}{dt}|_{t=0})$ and therefore
$ev_0(\frac{d\psi_t}{dt}|_{t=0}) \in \tau'(f)+\tau'(g)$, which contradicts
our
hypotheses.\\

\ni Since $\frac{dh_\lambda}{d\lambda}|_{\lambda=0} \not\in T{\cal {A}}_e h$
and $h$
has ${\cal {A}}_e$-codimension 1, $H$ is a versal unfolding of $h$.
\eop

\begin{corl}\label{sub}
If $h:(\CC^n,S)\to(\CC^p,0)$ is a multi-germ of ${\cal {A}}_e$-codimension 1,
then for every proper subset $S'$ of $S$, the restriction of $h$ to a
multi-germ $(\CC^n,S')\to(\CC^p,0)$ is stable.
\end{corl}
\Proof Let $S=S'\cup S''$ with $S'\cap S''=\emptyset$. Let $h'$ and $h''$
be the
multi-germs of $h$ at $S'$ and $S''$ respectively. Suppose that one of $h'$
and
$h''$ is not stable, say $h'$.
Then $h'$ has ${\cal
{A}}_e$-codimension 1. Since it is therefore
not a prism, by Lemma~\ref{tauprism}
$\tau'(h')=0$. As $h''$ is not a submersion, we may
choose
$v \in T_0 \Bbb {C}^p \backslash \tau'(h'')$. Extend $v$ to a vector field on
$\Bbb {C}^p$ and
integrate it to give a germ of a 1-parameter family $\psi_t$ of
diffeomorphisms of $(\Bbb {C}^p,0)$ satisfying the conditions of
Proposition \ref{tauspace}. Therefore $H$, as described in
Proposition \ref{tauspace}, is a versal unfolding of $h$. But then $id_\Bbb
{C}\times h'$ is a versal
unfolding of $h'$ and so $h'$ is stable, a contradiction.
Therefore $h'$ and $h''$ are stable.
\eop

A finite set
$E_1, \ldots, E_s$ of vector subspaces of a finite dimensional vector
space $F$ has
{\bf almost regular intersection} (with respect to $F$) if
$${\rm codim}(E_1 \cap \ldots \cap E_s)=
{\rm codim}E_1 + \cdots + {\rm codim}E_s
-1$$

\begin{lema}

$E_1, \ldots, E_s$ have almost regular intersection if and
only if the cokernel of
the natural mapping
$$F \to (F/E_1) \oplus \ldots \oplus (F/E_s)$$
has dimension 1.
\eop
\end{lema}

\begin{prop}\label{alreg}

Let $h=\{f,g\}$ be an ${\cal A}_e$-codimension 
1 multi-germ.
Then $\tau (f)$ and $\tau (g)$ have
almost regular intersection with respect to $T_0 \Bbb {C}^p$.

\end{prop}
\Proof Let $H$ be a versal unfolding of $h$. $H$ restricts to a versal
unfolding $F$
of $f$ and a versal unfolding $G$ of $g$. Since $f$ is stable,
$F$ is equivalent
to a prism on $f$ and hence
$$ T_0 \Bbb{C}^p / \tau(f) \cong T_0 (\Bbb{C} \times \Bbb{C}^p) / \tau(F)$$

\ni We have the following commutative diagram
$$\begin{array}{ccc}
T_0 \Bbb {C}^p & \rightarrow & \frac{T_0 \Bbb {C}^p}{\tau (f)}
\oplus \frac{T_0
\Bbb {C}^p}{\tau (g)}\\
\downarrow & & \downarrow\\
T_0 (\Bbb {C} \times \Bbb {C}^p) &
\rightarrow & \frac{T_0 (\Bbb {C} \times \Bbb
{C}^p)}{\tau (F)} \oplus \frac{T_0 (\Bbb {C} \times \Bbb {C}^p)} {\tau (G)}
\end{array}$$
in which the right hand map is bijective and the bottom map is
surjective by \ref{regular}.
So the top map has cokernel of codimension at most 1. Were it surjective,
then $\tau(f)$ and $\tau(g)$ would be transverse, and $h$ would be stable.
Hence the dimension of the cokernel is 1, proving the proposition.\eop

\begin{corl}\label{almostreg}

If $h$ is a multi-germ of ${\cal {A}}_e$-codimension 1 with branches $h^{(1)},
\ldots, h^{(r)}$, $r \geq 2$,
then $\tau (h^{(1)}),\ldots, \tau (h^{(r)})$ have almost regular
intersection with
respect to $T_0 \Bbb {C}^p$. \eop
\end{corl}

\begin{corl}\label{codtau}

Let $h=\{f,g\}$ have ${\cal A}_e$-codimension 1.
Then the codimension of $\tau
(f)
+\tau (g)$ in $T_0 \Bbb {C}^p$ is 1. \eop
\end{corl}

\ni It is natural to ask how we can tell when our codimension 1 multi-germ is
primitive.

\begin{prop}\label{primiaug}

Let $h=\{f,g\}$ be an ${\cal {A}}_e$-codimension 1 multi-germ,
and let $k=\dim_{\CC}$ $\tau(f)\cap\tau(g)$.
Then $h$ is a $k$-fold
augmentation of a primitive map-germ.
\end{prop}
\Proof By
Corollary~\ref{almostreg} we can choose $v \in T_0 \Bbb {C}^p \setminus
(\tau (f) + \tau (g))$. Choose a germ of a one parameter family $\psi_t$ of
diffeomorphisms of $(\Bbb {C}^p,0)$ such that
$ev_0(\frac{d\psi_t}{dt}|_{t=0})=v$.
Then choose a versal unfolding $H$ of $h$ as in
Proposition~\ref{tauspace}. If
$\Lambda$ is the first coordinate in the target
$\Bbb{C} \times \Bbb{C}^p$ of $H$
then $\tau(F)=\Bbb{C}\frac{\p}{\p\Lambda}
\oplus\tau(f)$ and $\tau(G)=
\Bbb{C}(\frac{\p}{\p\Lambda}+v)\oplus\tau(g)$. Since
$\tau (H)=\tau (F) \cap \tau (G)$ it follows that $\tau (H)=\tau (f) \cap
\tau
(g)$. Therefore, by Proposition~\ref{anastratum}, $H$ is a prism and by
Theorem \ref{augment} $h$ is an augmentation. \eop
\begin{corl}\label{decomp} Suppose that $h=\{f,g\}$ is a primitive
${\cal A}_e$-codimension 1 multi-germ. Then there is a decomposition
$$T_0 \Bbb {C}^p=\tau (f) \oplus \tau (g) \oplus \Bbb {C}v.$$
\end{corl}
\Proof Immediate from Corollary \ref{codtau} and 
Proposition \ref{primiaug}.\eop
\begin{example}\label{classim}
{\em Using \ref{primiaug} we classify codimension 1
multi-germs of immersions. If $f:\CC^n,S\to\CC^{n+1}$ has all of its
$r$ branches immersions, then the same is true of a 1-parameter 
versal unfolding $F$. As $F$ is stable, these $r$ branches meet in general
position, with intersection $L$ of dimension $n+1-r$. As $L=\tau(F)$, 
by \ref{primiaug} $F$ is the $n+1-r$-fold augmentation of a germ
$f_0:\CC^{r-1},S\to\CC^r,0$, evidently also consisting of $r$
immersions. As $f_0$ has ${\cal A}_e$-codimension 1, a little thought shows
that each $r-1$-tuple of its immersions is in general position (but see also
\ref{alreg}). It follows that $f_0$ is equivalent to the germ consisting of
a parametrisation of the $r$ coordinate hyperplanes, together with one extra
immersive branch $(x_1,\ldots,x_{r-1})\mapsto (x_1,\ldots,x_{r-1},\sum_ix_i)$.
This has a versal unfolding in which only the last immersion is deformed, to
$(x_1,\ldots,x_{r-1})\mapsto (x_1+t,\ldots,x_{r-1}+t,\sum_ix_i+t)$. Thus
$f$ is equivalent to the germ consisting of a parametrisation of the
first $r-1$ hyperplanes together with an additional immersion of the form
$$(x_1,\ldots,x_{r-1},u_1,\ldots,u_{n+r-1})\mapsto
(x_1+\sum_ju_j^2,\ldots,x_{r-1}+\sum_ju_j^2, \sum_ix_i+\sum_ju_j^2,u_1,\ldots,
u_{n-r+1}).$$
In the real case, the only change in the classification is
that $\sum_ju_j^2$ must be replaced by $\sum_j\pm u_j^2$, giving 
$(n-r+1)/2$ different classes if $n+r-1$ is even, or
$(n+r)/2$ if $n-r+1$ is odd.\\

\ni The second germ in the list shown in the right-hand column
in Figure 1 on page 10 is of this
type.}
\end{example}

In view of \ref{decomp}, by a change of coordinates we can arrange that 
the analytic stratum of $f$ becomes $\Bbb {C}^a \times \{0\}
\times \{0\}$,
that of $g$ becomes $\{0\} \times \Bbb {C}^b \times \{0\}$ and $v$
becomes $(0,0,1) \in \Bbb {C}^a \times \Bbb {C}^b \times \Bbb {C}$. We shall
suppose for the
remainder of this section that this change of coordinates has been
made.

We say that a multi-germ $f$ is transverse to a vector subspace
$V$ of $T_0 \Bbb
{C}^p$ if every branch of $f$ is transverse to $V$.
Our analysis of multi-germs $h=\{f,g\}$ from now on falls
into two cases, characterised by whether
$g$ is or is not transverse to $\tau (f)$.\\

\ni{\bf Case 1:}
$g$ is not transverse to $\tau (f)$.

\begin{lema}\label{morse}

A stable map germ of rank zero is
either a Morse singularity, or either
the domain or the codomain has dimension zero.
\eop
\end{lema}

\begin{prop}\label{morpri}
Let $h=\{f,g\}$ be a primitive ${\cal {A}}_e$-codimension 1 multi-germ,
and suppose that
$g$ is not transverse to $\tau (f)$.
Then 
\be
\item
if moreover $g$ and $f$ are transverse, it follows that
\be
\item
$g$ has precisely one branch, which
is either a prism on a Morse singularity or an immersion. 
\item
After a change of coordinates, $h$ takes the form
$$\cases{ f:(\CC^{n-1}\times\CC,S_0\times\{0\})\to(\CC^{p-1}\times\CC,0),
\ \ f(x,u)=(f_u(x),u)\cr
g:(\CC^{p-1}\times\CC^k,0)\to(\CC^{p-1}\times\CC,0),\ \ g(\lambda,v)=
(\lambda,\sum_j v_j^2)}
$$
where $f$ is an ${\cal A}_e$-versal unfolding of $f_0$; thus
$h=C_k(f_0)$.
In particular, $f\tr\tau(g)$.
\ee
\item
if $g$ and $f$ are not transverse, then $p=1$, and $f$ and $g$
are both Morse functions.
\ee
\end{prop}
\Proof
If $g$ has more than one branch, then by \ref{sub} the multi-germ
consisting of
$f$ together with any one branch $g^{(i)}$ of $g$ is stable. Hence
$\tau(g^{(i)})\tr\tau(f)$, so $g^{(i)}\tr\tau(f)$, so $g\tr\tau(f)$. This
contradiction implies that $g$ has only one branch.

Now suppose that $\mbox{Image}(dg(0))$ is bigger than $\tau(g)$. 
Then we can construct a 
1-parameter deformation $h_t$ of $h$ by fixing $f$ and composing $g$
with a 1-parameter rotation about $\tau(g)$, in such a way that for $t\neq 0$,
$g$ becomes transverse to $\tau(f)$. Since $\tau(g)$ remains
non-transverse to $\tau(f)$, $h_t$ is not stable even for $t\neq 0$. But
neither is it equivalent to $h=h_0$. This is impossible, since $h$ has
$\AA$-codimension 1.
Hence $\mbox{Image}(dg(0))=\tau(g)$, and so $g$ is a prism on a germ
of rank 0. By \ref{morse}, $g$ is either a prism on a Morse function
or an immersion.\\
\ni (2) The codimension of $\tau(g)$ is now $1$, so by Corollary \ref{decomp}
we must have $\tau(f)=\{0\}$. 
Thus, we have a decomposition of the target as
$\Bbb {C}^{p-1} \times \Bbb {C}$ where $\Bbb
{C}^{p-1} \times \{0\}$ is the analytic stratum of $g$.
There is a neighbourhood
$U$ of 0 in $\Bbb {C}^{p-1}$ such that for all $u \in U$, the pullback of $g$
along the inclusion of the
subset $\{u\} \times \Bbb {C}$ is a Morse singularity
and so by a coordinate change in the source
we can
reduce this pullback to the form $\sum^m_{i=1} x_i^2$. In fact 
the changes of coordinates in the source
depend analytically on $u$ and
so together they give a change of coordinates in the
source which reduces $g$ to the
form
$$\begin{array}{rcl}
\Bbb {C}^{p-1} \times \Bbb {C}^k & \to & \Bbb {C}^{p-1} \times \Bbb {C}\\
(\lambda,v_1,\ldots,v_k) & \mapsto & (\lambda,\sum^m_{j=1}v_j^2) \end{array}$$

Now suppose that $f$ is transverse to $g$. Then 
by a change of coordinates in the source of 
$f$ we can now bring $f$ to the desired form. 
Evidently $f$ is now a stable 1-parameter unfolding of $f_0$, so
we can view $h$ as $C_k(f_0)$; finally, by Theorem \ref{conc1} 
$${\cal A}_e-\mbox{codim}(f_0)={\cal A}_e-\mbox{codim}(h)=1.$$

On the other hand, if $f$ is not transverse to $g$ then we can apply the
previous argument with the roles of $f$ and $g$ reversed, to conclude
that $p=1$ and thus that $f$ and $g$ are both Morse singularities.
\eop
\begin{example} {\em The germ $f_0$ of Example \ref{classim} is obtained
(up to $\cal A$-equivalence) by
applying the concatenation operation $C_0$ (defined 
using Theorem \ref{indep2}) $r-1$ times to the bi-germ
consisting of coincident embeddings of two copies of
$\CC^0$ in $\CC$.}
\end{example}

To complete our analysis of codimension 1 multi-germs, by \ref{morpri} 
it remains to consider only\\ 
{\bf Case 2: $f\tr\tau(g)$ and $g\tr\tau(f)$}.
Recall that we were able to
decompose the target $\CC^p$ as $\CC^a\times\CC^b\times\CC$, with
$\tau(f)=\CC^a\times\{0\}\times\{0\}$ and 
$\tau(g)=\{0\}\times\CC^b\times\{0\}$. Let $z_1,\ldots,z_{a+b+1}$ be
coordinates on $\CC^a\times\CC^b\times\CC$. 
Since $f$ is transverse to $\tau(g)$, we can take $z_{a+b+1}\circ f$ as 
a coordinate, $u$, on the domain of $f$, and similarly, as $g$ is transverse to
$\tau(f)$, we can take $v=z_{a+b+1}\circ g$ as a coordinate on the
domain of $g$. A coordinate change now brings $\{f,g\}$ to the form
$$\cases{(X,y,u)\mapsto (X,f_u(y),u)\cr
(x,Y,u)\mapsto (g_{u,Y}(x),Y,u)\cr}.
$$
Note that we have reduced $f$ to a prism on a 1-parameter unfolding
(using the fact that $\CC^a\times\{0\}\times\{0\}$ is the analytic stratum of
$f$), but that we have not done the same for $g$ --- yet. A naive coordinate
change to reduce $g$ to a prism on a 1-parameter unfolding would take $f$
out of its normal form. Nevertheless, we claim that $h$ is $\cal A$-equivalent
to a binary concatenation of two ${\cal A}_e$-codimension 1 germs, as 
described in Section \ref{concatenation}. As a first step, we prove:
\begin{lema}\label{gebar}
Suppose that $h$ is an ${\cal A}_e$-codimension 1 germ in the
semi-normal form
$$\cases{(X,y,u)\mapsto (X,f_u(y),u)\cr
(x,Y,u)\mapsto (g_{u,Y}(x),Y,u)\cr}.
$$
Then 
\be
\item
the ${\cal A}_e$-codimension of the germs $g_0$ and $f_0$ 
is equal to $1$, 
and the germs $\overline{g}:(x,v)\mapsto (g_{v,0}(x),v)$ and
$\overline{f}:(y,u)\mapsto (f_u(y),u)$ are ${\cal A}_e$-
versal unfoldings of $g_0$ and $f_0$. 
\item
If also $h$ is primitive, then so are $g_0$ and $f_0$.
\ee
\end{lema}
\Proof We give the proof for $g_0$ and $\overline{g}$; the proof of $f_0$ and
$\overline{f}$ is identical.\\
{\bf Step 1:} The unfolding $H$ of $h$ given by
$$\cases{(X,y,u,v)\ \ \stackrel{F}{\mapsto}\ \ (X,f_u(y),u+v,v)\cr
(x,Y,s,v)\ \ \stackrel{g\times\mbox{id}_{\CC}}{\mapsto}\ \ (g(x,Y,s),v)}$$
is ${\cal A}_e$-versal. For it is not infinitesimally trivial, and $h$
has ${\cal A}_e$-codimension 1.\\
{\bf Step 2:}
Let $G$ be an unfolding of $g_0$, and let 
$\tilde{G}$ be the direct sum unfolding of $G$ and $g$. 
Clearly $G$ can be induced from $\tilde{G}$.
Consider the unfolding $\tilde{H}$ of $h$, given by $\tilde{H}=\{F\times
\mbox{id}_{\CC^d},\tilde{G}\}$. As the 1-parameter unfolding
$H$ of $h$ is versal, $\tilde{H}$ must be isomorphic to an unfolding
induced from $H$. This means $\tilde{G}$ is isomorphic to an unfolding
induced from $g\times\mbox{id}_{\CC}$. Any such unfolding is isomorphic
to an unfolding induced from $g$. Hence $g$ is a versal unfolding of
$g_0$. The Kodaira-Spencer map of $g$, 
from  
$T_0\CC^b\times\CC$
to the ${\cal A}_e$-normal space of $g_0$, is therefore surjective.
But as $g$ is trivial along $\{0\}\times\CC^b\times\{0\}$, the Kodaira-
Spencer map is identically zero along $\CC^b\times\{0\}$. Hence the restriction
of the Kodaira-Spencer map to $\{0\}\times\CC$ is surjective, and 
$\overline{g}$ is ${\cal A}_e$-versal.\\

\ni 
If also $h$ is primitive, then $\tau(f)\cap\tau(g)=\{0\}$, by \ref{primiaug}, 
and
so the analytic stratum of the versal unfolding $(x,v)\mapsto (x,g_v(x))$
must be reduced to $\{0\}$ also. It follows from \ref{augment}
that $g_0$ must be primitive.
\eop
\begin{corl}\label{banana}
Suppose that $\{f,g\}$ is a multi-germ of ${\cal A}_e$-
codimension 1,  with $f$ transverse to $\tau(g)$ and $g$ transverse to
$\tau(f)$. Then the pull-back of $f$ by $\tau(g)$, and the pullback of
$g$ by $\tau(f)$, are both germs of ${\cal A}_e$-codimension 1. 
\end{corl}
\Proof When $\{f,g\}$ is put in the semi-normal form of the Proposition,
these pull-backs are just $f_0$ and $g_0$, and the proposition establishes
that they have ${\cal A}_e$-codimension 1. However, the statement is
evidently independent of choice of coordinates.
\eop
We would like to be able to put the germ $h=\{f,g\}$ of \ref{gebar}
into a normal 
form, 
$$\cases{(f:(X,y,u)\mapsto (X,f_u(y),u)\cr
G:(x,Y,v)\mapsto ((g_v(x),Y,v)}
;$$
but it is not clear that this is always possible. The problem is as follows:
now that we have established
that $(x,u)\mapsto (g_{0,u}(x),u)$ is a versal unfoding of $g_0$, it follows
that there exists a submersion $\gamma:\CC^b\times\CC\to\CC$, and germs of
families of diffeomorphisms $\phi_{Y,u},\psi_{Y,u}$ such that
$$g_{Y,u}=\psi_{Y,u}\circ\g_{0,\gamma(Y,u)}\circ \phi_{Y,u};$$
nevertheless, in order to transform $h$ from its semi-normal form to the
desired normal form, the $\psi_{Y,u}$ and $\phi_{Y,u}$ would have to satisfy
the stronger requirement that
$$g_{Y,u}=\psi_{Y,u}\circ\g_{0,u}\circ \phi_{Y,u}.$$
This can be done under certain assumptions of quasihomogeneity, which we now
explain.\\

\ni A map $f:\Bbb {C}^n \to \Bbb {C}^p$ is {\bf weighted homogeneous}
if there
are positive integers $\omega_1, \ldots, \omega_n$ (the weights) and $d_1,
\ldots,d_p$ (the degrees) such that for $\mu \in \Bbb {C}$,
$f(\mu^{\omega_1}x_1,\ldots,\mu^{\omega_n}x_n)=(\mu^{d_1}f_1(x),\ldots,
\mu^{d_p}f_p(x))$. A germ $f$ is {\bf quasihomogeneous} if it is 
$\cal{A}$-equivalent to a weighted homogeneous map-germ. A multi-germ is
quasihomogeneous if its branches are quasihomogeneous with the same degrees.

Let $f\!:\!(\Bbb {C}^n,S) \to (\Bbb {C}^p,0)$\, be a quasihomogeneous
multi-germ of
${\cal {A}}_e$-codimension 1. When $(n,p)$ are in the range of nice
dimensions, we can find a
quasihomogeneous versal unfolding $F(\lambda,x)=(\lambda,f_\lambda(x))$ of
$f$ such
that the degree $r$ of the unfolding parameter is positive. In fact, if the
degree is
non-positive, then $F$ is topologically trivial and therefore $f$ is
topologically
stable. But this is a contradition since in the nice dimensions topological
stability is equivalent to stability. Let $r,d_1,\ldots,d_p$ be 
the degrees of the components 
of $F$
and let $r, w_1^{(i)}, \ldots, w_{n(i)}^{(i)}$ be the weights in the source of
the
$i^{th}$ branch of $F$. For $\mu \in \Bbb {C}$ define $\psi_\mu:\Bbb {C}^p
\to
\Bbb {C}^p$ by $\psi_\mu(y_1,\ldots,y_p)=(\mu^{d_1}y_1,\ldots,
\mu^{d_p}\,y_p)$ and define $\Psi_\mu:\Bbb {C}^{p+1} \to \Bbb {C}^{p+1}$
by
$\Psi_\mu(\lambda,y)=(\mu^r \lambda,\psi_\mu(y))$. Let
$\phi_\mu^{(i)}$ and
$\Phi_\mu^{(i)}$ be the analogues of these maps in the source of the $i^{th}$
branch of $f$ and $F$ respectively. If $\phi_\mu$ has branches
$\phi_\mu^{(i)}$
then
$$f_{\mu^r \lambda} \circ \phi_\mu=\psi_\mu \circ f_\lambda$$

\begin{lema}\label{change}

Let $\bar{f}$ and $\tilde{f}$ be quasihomogeneous $\cal {A}$-equivalent 
multi-germs from $\Bbb
{C}^n$ to
$\Bbb {C}^p$ ($(n,p)$ nice dimensions) of ${\cal {A}}_e$-codimension 1. 
Let $\bar{F}(\l_1,\ldots,\l_d,x)\!=\!(\l_1,\ldots,\l_d,
\bar{f}_{\l_1,\ldots,\l_d}(x))$ be a
versal unfolding of $\bar{f}$ with analytic stratum $\{0\} \times \Bbb
{C}^{d-1}
\times \{0\}$ and let $\tilde{F}(\mu,x)=(\mu,\tilde{f}_\mu(x))$ be a
versal
unfolding of $\tilde{f}$. Then there are families of diffeomorphisms
$\alpha_\lambda$ of $\Bbb {C}^n$ and $\beta_\lambda$ of $\Bbb {C}^p$,
$\lambda \in
\Bbb {C}^d$, such that the following diagram commutes
$$\begin{array}{ccc}
\Bbb {C}\times \Bbb {C}^{d-1}\times \Bbb {C}^n & \stackrel{\bar{F}}
{\rightarrow}
& \Bbb {C}\times \Bbb {C}^{d-1}\times \Bbb {C}^p\\ \,\,\,\downarrow \alpha
& &
\,\,\,\downarrow \beta\\ \Bbb {C} \times \Bbb {C}^n &
\stackrel{\tilde{F}}{\rightarrow} & \Bbb {C}\times \Bbb {C}^p
\end{array}$$
where $\alpha(\mu,\nu,x)=(\mu,\alpha_{(\mu,\nu)}(x))$ and
$\beta(\mu,\nu,y)=(\mu,\beta_{(\mu,\nu)}(y))$.

\end{lema}
\Proof We may suppose that $\tilde {f}$ and $\tilde{F}$ are
quasihomogeneous as
maps.

Let $\phi$ and $\psi$ be diffeomorphisms such that $\psi \circ \bar{f}=
\tilde{f}
\circ \phi$. Then $F'=(id_{\Bbb {C}^d} \times \psi) \circ \bar{F} \circ
(id_{\Bbb
{C}^d} \times \phi)^{-1}: \Bbb {C}^d \times \Bbb {C}^n \to \Bbb {C}^d
\times \Bbb
{C}^p$ is a versal unfolding of $\tilde{f}$ with analytic stratum
$\{0\}\times
\Bbb {C}^d \times \{0\}$.

Since $\tilde{F}$ is a miniversal unfolding, there is a submersion
$\gamma: \Bbb
{C}^d \to \Bbb {C}$ and there are families of diffeomorphisms
$\bar{\phi}_\lambda$
of $\Bbb {C}^n$ and $\bar{\psi}_\lambda$ of $\Bbb {C}^p$, $\lambda \in \Bbb
{C}^d$, such that the following diagram commutes.
$$\begin{array}{ccc}
\Bbb {C}\times \Bbb {C}^{d-1}\times \Bbb {C}^n & \stackrel{F'}
{\rightarrow} &
\Bbb {C}\times \Bbb {C}^{d-1}\times \Bbb {C}^p\\ \,\,\,\,\,\,\,\,\,\downarrow
\Gamma\times \bar{\phi}_\lambda & &\,\,\,\,\,\,\,\,\, \downarrow \Gamma
\times
\bar{\psi}_\lambda\\ \Bbb {C} \times \Bbb {C}^{d-1}\times \Bbb {C}^n &
\stackrel
{\tilde{F} \times id_{\Bbb {C}^{d-1}}}{\rightarrow} & \Bbb {C}\times \Bbb
{C}^{d-1} \times \Bbb {C}^p \end{array}$$
where $\Gamma(\mu,\nu)=(\gamma(\mu,\nu),\nu)$.

We have $\gamma^{-1}(0)=\{0\}\times \Bbb {C}^{d-1}$, so $\Gamma$ is
a diffeomorphism
by the inverse function theorem. Since $\Gamma$ commutes with projection onto
$\Bbb {C}^{d-1}$, $\Gamma^{-1}$ does also, so there exists $\gamma',
\gamma'':
\Bbb {C}^d \to \Bbb {C}$ such that
$\Gamma^{-1}(\mu,\nu)=(\gamma'(\mu,\nu),\nu)$
and $\gamma'=\mu \gamma''$ where
$\mu:\Bbb {C} \times \Bbb {C}^{d-1} \to \Bbb {C}$ is the projection onto
the first
coordinate. Also $\gamma'$ is a submersion and so $\gamma''$ is non-zero in a
neighbourhood of the origin. We have
$$(\gamma' \times \psi_{\sqrt[r]{\gamma''}})\circ (\tilde{F} \times
id_{\Bbb
{C}^{d-1}})=\tilde{F} \circ (\gamma' \times \phi_{\sqrt[r]{\gamma''}})$$
where $\phi_{\sqrt[r]{\gamma''}}$ and
$\psi_{\sqrt[r]{\gamma''}}$ are as stated just before this
proposition. Thus
the following diagram commutes
$$\begin{array}{ccc}
\Bbb {C}\times \Bbb {C}^{d-1}\times \Bbb {C}^n & \stackrel{F'}
{\rightarrow} &
\Bbb {C}\times \Bbb {C}^{d-1}\times \Bbb {C}^p\\ \,\,\,\downarrow
\bar{\alpha} &
&\,\,\, \downarrow \bar{\beta}\\ \Bbb {C} \times \Bbb {C}^n &
\stackrel{\tilde{F}}{\rightarrow} & \Bbb {C}\times \Bbb {C}^p
\end{array}$$
where $\bar{\alpha}=(\gamma' \times \phi_{\sqrt[r]{\gamma''}})
\circ (\Gamma \times \bar{\phi}_\lambda)$ and $ \bar{\beta}=((\gamma' \times
\psi_{\sqrt[r]{\gamma''}}) \circ (\Gamma \times \bar{\psi}_\lambda)$.

Now the proposition follows by choosing $\alpha=\bar{\alpha} \circ (id \times
\phi)$ and $\beta=\bar{\beta}\circ (id \times \psi)$. \eop

Now we can continue with the task of 
reducing a primitive 
${\cal A}_e$-codimension 1 germ in the semi-normal
form
$$\cases{(X,y,u)\mapsto (X,f_u(y),u)\cr
(x,Y,v)\mapsto (g_{Y,v}(x),Y,v)}
$$ 
to the normal form
$$\cases{(X,y,u)\mapsto (X,f_u(y),u)\cr
(x,Y,v)\mapsto (g_v(x),Y,v)}.
$$
We make the additional hypothesis that $g_0$ is quasihomogeneous, and 
is not topologically stable. Then in appropriate coordinates 
it has an ${\cal A}_e$-versal unfolding
whose unfolding parameter has positive weight. Thus we can apply 
\ref{change}, to deduce that the unfolding $g$ of $g_0$ is isomorphic to
a prism on the unfolding $\overline{g}:\ (x,u)\mapsto (g_{0,u}(x),u)$.
That is, there are diffeomorphisms $\Phi:(\CC^\alpha\times\CC^b\times\CC,T)
\to (\CC^\alpha\times\CC^b\times\CC,T)$, 
of the form $\Phi(x,Y,u)=(\phi_{Y,u}(x),Y,u)$, 
and $\Psi:(\CC^a\times\CC^b\times\CC,0)\to 
(\CC^a\times\CC^b\times\CC,0)$ of the form $\Psi(X,Y,u)=(\psi_{Y,u}(X),Y,u)$,
such that $g_{Y,u}(x)=\psi_{Y,u}\circ g_{0,u}\circ \phi_{Y,u}$.
Composing with $\Phi$ in the source of $g$, and $\Psi$ on the
target of $\{f,g\}$, we bring $\{f,g\}$ to the form
$$\cases
{(X,y,u)\mapsto (\psi_{f_u(y),u}(X),f_u(y),u)\cr
(x,Y,v)\mapsto (g_{0,v}(x),Y,v)}
$$
and now if we take the first $a$ coordinates of $\Psi\circ f$ as new
coordinates on the domain of $f$, we bring $\{f,g\}$ to the desired normal
form.
We have proved
\begin{teo}\label{bin}
If $h=\{f,g\}$ is a multi-germ of ${\cal A}_e$ codimension
1, in which $f$ is transverse to $\tau(g)$ and $g$ is transverse to $\tau(f)$,
and if either the pullback of $f$ by $\tau(g)$ or the pullback of $g$ by 
$\tau(f)$ is quasihomogeneous and not topologically stable, then $\{f,g\}$
is equivalent to a binary concatenation $B(f_0,g_0)$; that is, to a germ of the form
$$\cases
{(X,y,u)\mapsto (X,f_u(y),u)\cr
(x,Y,v)\mapsto (g_v(x),Y,v)}.
$$
\eop
\end{teo}
We now summarise the results of this section:\\

\begin{teo}\label{classi} 
Let $h=\{f,g\}$ be a primitive ${\cal A}_e$-codimension 1 map-germ in the
nice dimensions (with no submersive branches).
Then $f$ and $g$ are both stable (\ref{sub}). 
\be
\item
If $f$ and $g$ are not transverse, then (\ref{morpri}) $h$ is equivalent to
$$
\cases{
(x_1,\ldots,x_n)\mapsto \sum_ix_i^2\cr
(y_1,\ldots,y_m)\mapsto \sum y_j^2
.}
$$
Now assume $f\tr g$.\\
\item
If $g$ is not transverse to $\tau(f)$, 
then (\ref{morpri}) $f$ is transverse to
$\tau(g)$, and $h$ is equivalent to
$$
\cases{
(x_1,\ldots,x_n,u)\mapsto (f_u(x),u)\cr
(\lambda_1,\ldots,\lambda_{p-1},v_1,\ldots,v_k)\mapsto (\lambda, \sum_iv_i^2)
}
$$
(so $\{f,g\}$ is equivalent to $C_k(f_0)$).
\item If $g\tr\tau(f)$  and $f\tr\tau(g)$,
then (\ref{bin}) $\{f,g\}$ is equivalent to a germ of the form
$$
\cases{
(X,y,u)\mapsto (X,f_u(y),u)\cr
(x,Y,v)\mapsto (g_{Y,v}(x),Y,v)
}
$$
where the target is decomposed as $\CC^a\times\CC^b\times\CC$, and 
$f_0$ and $g_0$ are primitive. If also
the pullback of $g$ by $\tau(f)$ or the pullback of $f$ by $\tau(g)$
is quasihomogeneous then $\{f,g\}$ is equivalent to
$$
\cases{
(X,y,u)\mapsto (X,f_u(y),u)\cr
(x,Y,v)\mapsto (g_v(x),Y,v),
}
$$
i.e. to $B(f_0,g_0)$. 
\ee
\end{teo}\eop
\begin{obs}\label{real}{\em If we replace $\Bbb {C}$ by $\Bbb {R}$ and
analytic
maps by smooth ones, then the results obtained so far still hold modulo the
following alterations:
in the real case we define two augmentations:
$ A_F^+(\lambda,x)=(\lambda,f_{\lambda^2}(x))$ and
$ A_F^-(\lambda,x)=(\lambda,f_{-\lambda^2}(x))$.
In the proof of Proposition~\ref{change}, if $\omega$ is
even then
we cannot necessarily define $\sqrt[\omega]{}$\, properly. Consequently we
may
have to define $\alpha(\mu,\nu,x)=(-\mu,\alpha_{(\mu,\nu)}(x))$ and
$\beta(\mu,\nu,y)=(-\mu,\beta_{(\mu,\nu)}(y))$ in order to get the diagram to
commute.}
%
\end{obs}
\section{Topology}\label{topology}

Let $f\!:\!(\Bbb {C}^n,S) \to (\Bbb {C}^p,0)$\, ($n\geq p-1$, $(n,p)$ nice
dimensions and $S$ a finite set) be a finitely $\cal {A}$-determined
multi-germ. A
{\bf stabilisation} of $f$ is a 1-parameter unfolding $F\!:\!(\Bbb {C}
\times \Bbb {C}^n, \{0\} \times S) \to (\Bbb {C} \times \Bbb {C}^p, (0,0))$\,
with the property that there is a representative $F\!:\!U \to V$\, (we
shall use
the same letter) and a
positive real number $\delta$ such that for $\lambda \in B_\delta (0)\setminus
\{0\}$,
the map
$f_\lambda:U_\lambda \to V_\lambda$ \, is infinitesimally stable (here
$U_\lambda=U \cap (\{ \lambda \} \times \Bbb {C}^n)$ and $V_\lambda=V \cap
(\{
\lambda \} \times \Bbb {C}^p)$), $F|_{\sum(F)}$ is proper, finite to one and
generically one to one, and that $F^{-1}(0,0)\cap \sum(F)=\{0\}\times S$. It
follows that the discriminant $D(F)$ of $F$ is a closed analytic subset of
$V$. The mapping $f_\lambda$ is a {\bf stable perturbation} of $f$.

Consider the canonical stratification of $D(F)$ and choose $\epsilon >0$
such that
for all $\epsilon'$ with $0< \epsilon' \leq \epsilon$, $D(f) \cong D(F) \cap
(\{0\} \times \Bbb {C}^p)$ is stratified transverse to the sphere
$S_{\epsilon'}
\subset \Bbb {C}^p$ of centre 0 and radius $\epsilon'$. Such $\epsilon$ is
called
a {\bf Milnor radius} for $D(f)$. By Thom's First Isotopy Lemma, $D(f) \cap
B_\epsilon$ is a cone on its boundary $D(f) \cap S_\epsilon$. It follows that
there is a $\delta >0$ such that for $\lambda \in B_\delta \subseteq \Bbb
{C}$,
$D(F)$ is stratified transverse to $\{\lambda\} \times S_\epsilon $ (we
call such
a $\delta$ a {\bf perturbation limit} for $F$ with respect to
$B_\epsilon$). For
$\lambda \in B_\delta$, the {\bf discriminant} of $f_\lambda$ is defined
to be
$D(f_\lambda) \cap B_\epsilon$, or, in other words, $D(F) \cap (\{ \lambda \}
\times B_\epsilon)$.

For $\epsilon_1, \ldots, \epsilon_p > 0$ define the set $P_{\epsilon_1,
\ldots,
\epsilon_p}(0)$ to be the polycylinder $\{(y_1, \ldots, y_p) \in \Bbb
{C}^p \,/\, |y_i|< \epsilon_i\,\, \forall i\}$. We shall also use the term
``Milnor
radius for
$D(f)$'' for an $\epsilon > 0$ such that for all $\epsilon_1, \ldots,
\epsilon_p$
with $0 < \epsilon_i < \epsilon$ \,($ \forall i$), $D(f)$ is stratified
transverse
to the boundary of the polycylinder $P_{\epsilon_1, \ldots,
\epsilon_p}(0)$. The
results described above apply with such a polycylinder replacing
$B_\epsilon $ and
the discriminant defined this way is the same.

Let $\pi:D(F) \to \Bbb {C}$\, be the projection to the parameter space
$\Bbb {C}$.
It follows by \cite{damon-mond} that $\pi$ induces a
locally
trivial fibration
$$((B_\delta \backslash \{0\}) \times B_\epsilon) \cap D(F) \to B_\delta
\backslash \{0\}$$

\begin{lema}\label{susp}

Let $A,B$ be contractible open subsets of a 
topological space $X$, and $A',B'$ be contractible
open subsets of $X'$. Suppose that
$A \cap B$ and $A' \cap B'$ are homotopy equivalent, and moreover that
$A \cap B$ has collared neighbourhoods in both $A$ and $B$, and $A' \cap
B'$ has
collared neighbourhoods in both $A'$ and $B'$. Then $A \cup B$ and $A'
\cup B'$
are homotopy equivalent.
\eop
\end{lema}

Suppose $f$ has ${\cal {A}}_e$-codimension 1. Let $F(\lambda,x)=(\lambda,
f_\lambda(x))$ be a proper representative of a miniversal unfolding of
$f$. For
$\mu \in \Bbb {C}$ define $g_\mu (\lambda,x) =(\lambda,
f_{\lambda^2+\mu}(x))$.
Then $G(\mu, \lambda,x)= (\mu,\lambda, f_{\lambda^2+\mu}(x))$ is a proper
representative of a miniversal unfolding of $g=A_F f$.

\begin{teo}\label{suspdis}

With the above notation, for $\mu\neq 0\neq \lambda$ 
the discriminant of $g_\mu$ is homotopy
equivalent to the
suspension of the discriminant of $f_\lambda$.

\end{teo}
\Proof Let $\epsilon >0$ be a Milnor radius for both $f$ and $F$,
also let
$\delta >0$ be a perturbation limit for $F$ with respect to $P_{\epsilon,
\ldots,
\epsilon}(0) \subseteq \Bbb {C}^p$.

Let $\epsilon' >0$ be a Milnor radius for $g$ and let $\delta' >0$
be a
perturbation limit for $G$ with respect to $P_{\epsilon'', \epsilon, \ldots,
\epsilon}(0) \subseteq \Bbb {C}^{p+1}$, where $\epsilon''=min\{ \epsilon,
\sqrt{\delta/2}\}$.

Fix $\mu_0 \in \Bbb {C}$ and consider
$$\pi:D(g_{\mu_0})\cap P_{\epsilon'', \epsilon, \ldots, \epsilon}(0) \to \Bbb
{C}$$
be the projection onto the first coordinate. For a convenient
choice of
$\mu_0$ we have

\vspace{.1in}

\noindent (i) the fibre of $\pi$ over $\lambda \in B_{\epsilon''}(0)$ is
naturally
homeomorphic to $D(f_{\lambda^2+\mu_0}) \cap P_{\epsilon, \ldots,
\epsilon}(0)$
which is the discriminant of $f_{\lambda^2+\mu_0}$.

\vspace{.1in}

\noindent (ii) Suppose that the square roots of $-\mu_0$ are in
$B_{\epsilon''}(0)$, say $a$ and $b$. Then the restriction of $\pi$ to
$\pi^{-1}(B_{\epsilon''}(0)\backslash \{a,b\})$ is a locally trivial
fibration.

\vspace{.1in}
Let $A$ and $B$ be contractible open subsets of $B_{\epsilon''}(0)$ with
contractible (non-empty) intersection such that $a \in A \backslash B$ and
$b \in
B \backslash A$. 
\vskip 10pt
\vskip 10pt
By standard arguments 
we can conclude that $\pi^{-1}(A \cup B)$ is homotopy
equivalent
to the discriminant of $g_{\mu_0}$ and $\pi^{-1}(A \cap B)$ is homotopy
equivalent
to the discriminant of $f_\lambda$; we 
can also assume that $\pi^{-1}(A \cap B)$ is collared in
both $\pi^{-1}(A)$ and $\pi^{-1}(B)$. 
Since the suspension of any space $D$ can be divided into two
contractible subspaces whose intersection has collared neighbourhoods and is
homotopy equivalent to $D$, by Lemma~\ref{susp} we have
only to
prove that $\pi^{-1}(A)$ and $\pi^{-1}(B)$ are contractible. At $a$,
$\gamma(\lambda)=\lambda^2+\mu_0$ is a 
diffeomorphism and induces a homeomorphism
between $\pi^{-1}(\gamma^{-1}(B_{\delta''}(0)))$ and $D(F) \cap
P_{\delta'',\epsilon, \ldots, \epsilon}(0)$ for some $\delta'' >0$. Therefore
$\pi^{-1}(A)$ is contractible since it is homeomorphic to
$\pi^{-1}(\gamma^{-1}(B_{\delta''}(0)))$ and $D(F) \cap P_{\delta'',\epsilon,
\ldots, \epsilon}(0)$ is a cone. Similarly, $\pi^{-1}(B)$ is contractible.
\eop

We now determine the homotopy-type of the discriminant of a stable
perturbation of a concatenation.

%
%

\begin{prop}\label{tipo} Let $f_0$ be a multi-germ of finite
${\cal A}_e$-codimension, which has a 1-parameter stable unfolding $F$.
The discriminant of a stable perturbation of the multi-germ
$C_k(f_0)$ (i.e. $\{F,g\}$, where $g(y,v)=(y,\sum v_i^2)$) is 
homotopy-equivalent to the suspension of the discriminant of
a stable perturbation of $f_0$.
\end{prop}

%
\Proof 
A stable perturbation $h_\lambda$ of $h$ has
branches $F$ and $g_\lambda(y,v)=(y,\sum v_i^2+\lambda)$. 
The discriminant of
$h_\lambda$ is
the union of two contractible spaces: the discriminant of $F$ and the
discriminant
of $g_\lambda$. The intersection of these sets is the discriminant of
$\tilde{f}_\mu$, which is a stable perturbation of $\tilde{f}_0$.
The proposition now follows from \ref{susp} in the same way 
as Theorem~\ref{suspdis}. 
\eop
In order to deal with the discriminant of a binary concatenation
$B(f_0,g_0)$, we need some
topological results.

Let $X$ and $Y$ be topological spaces. The {\bf join} of $X$ and $Y$, $X
\ast Y$,
is the space $(X \times Y \times I)/\sim$\, where $(x,y,\lambda) \sim
(x',y',\lambda')$ if and only if either $\lambda=\lambda'=0$ and $y=y'$ or
$\lambda=\lambda'=1$ and $x=x'$.

\begin{lema}\label{join}

If $X_1$ is homotopy equivalent to $X_2$ and $Y_1$ is homotopy equivalent
to $Y_2$
then $X_1 \ast Y_1$ is homotopy equivalent to $X_2 \ast Y_2$.
\eop
\end{lema}

\begin{corl}\label{susphom}

If $X_1$ is homotopy equivalent to $X_2$ then $S(X_1)$ is homotopy
equivalent to
$S(X_2)$.
\eop
\end{corl}

\begin{prop}\label{tipoiii}

Suppose that $h=B(f_0,g_0)$ is a binary concatenation, 
$$\cases{
(X,y,u)\ \stackrel{f}{\mapsto}\ (X,f_u(y),u)
\cr
(x,Y,v)\ \stackrel{g}{\mapsto}\ (g_v(x),Y,v)}
$$ of germs $f_0$ and $g_0$ of finite codimension, as described 
in Theorem \ref{conc2}. 
Let $H$ be the stable unfolding of $h$ given by
$$\cases{
(X,y,u,t)\ \stackrel{F}{\mapsto}\ (X,f_u(y),u+t,t)
\cr
(x,Y,v,t)\ \stackrel{g\times\mbox{id}_{\CC}}{\mapsto}\ (g_v(x),Y,v,t)}
.$$
Then for $t\neq 0$ 
the discriminant of the stable perturbation $h_t$ of $h$ is homotopy
equivalent to the suspension of $D(\tilde {f}_{-t}) \ast D(\tilde {g_t})$,
and thus 
$\mu_\Delta(h)=\mu_\Delta(\tilde {f_0})\times\mu_\Delta(\tilde {g_0})$.
\end{prop}
\Proof 
The discriminant of $h_t$ is the union of the (contractible)
discriminants of 
$(X,y,u)\mapsto (X,f_u(y),u+t)$ and $(X,y,v)\mapsto (g_v(x),Y,v)$.
It is prefereble to re-parametrise the first as the image of
$(X,y,u)\mapsto (X,f_{u-t}(y),u)$. 
Call these two spaces $D_1$ and $D_2$. By \ref{susp}, $D_1\cup D_2$ is
homotopy-equivalent to the suspension of 
$D_1\cap D_2$. Let $\epsilon>0$ be a Milnor radius for $f_0$ and $g_0$,
and let $P_f=P_{\epsilon,\ldots,\epsilon}(0)\subset \CC^b$
and $P_g=P_{\epsilon,\ldots,\epsilon}(0)\subseteq\CC^a$. 
Thus, we have to show that inside
a suitable Milnor polycyclinder $P_f\times P_g\times B(0,\epsilon')\subset
\CC^a\times\CC^b\times\CC$, and for $0 < |t| < \delta << \epsilon$,
$D_1\cap D_2$ is homotopy-equivalent
to the join of $D(f_{-t})\cap P_f$ and $D(g_t)\cap P_g$. 
This follows by a standard argument from the following three facts:
\begin{enumerate}
\item
The projection $\pi_h:\CC^a \times \CC^b \times \CC \to
\CC$\,
induces a locally trivial fibration
$$D_1\cap D_2 \cap \pi_h^{-1}(B_{\delta} \backslash \{0,t\}) \to B_{\delta}
\backslash \{0,t\}$$
whose fibre is homotopy equivalent to $D(g_t)\times D(f_{-t})$.
\item
The fibre of $\pi_h$ over $t$ is
$D(g_t)\times D(f_0)$;
because $D(f_0)$ is contractible, this
is homotopy-equivalent to
$D(g_t)$.
\item
The fibre of $\pi_h$ over $0$ is
$D(g_0)\times D(f_{-t})$;
because $D(g_0)$ is contractible, this
is homotopy-equivalent to
$D(f_{-t})$.
\end{enumerate}

Let $[0,t]$ denote the line-segment joining $0$ and $t$ in $\CC$.
Clearly there is a
deformation-retraction
$B_{\delta}\ \to\ [0,t]$; since $\pi_h$ is locally trivial on the
complement of
$[0,t]$, this lifts to a deformation-retraction
$D=\pi_h^{-1}(B_{\delta}) \ \to\ \pi^{-1}([0,t]$. By (1),(2),(3) above,
$\pi^{-1}([0,t]$ is homotopy-equivalent to
$D(g_t)\ast D(f_{-t})$.
 \eop

Now we consider the real case. Let $f\!:\!(\RR^n,S) \to (\RR^p,0)$ ($n
\geq p-1$, $(n,p)$ nice dimensions) be a multi-germ of
${\cal{A}}_e$-codimension 1
and let $F(\l,x)=(\l,f_\lambda(x))$
be a miniversal unfolding. Up to homeomorphism, there are two
(possibly equivalent) choices for the discriminant of $f_\lambda$: one
with positive $\lambda$ and one with negative $\lambda$. We shall call these
$D^+(f)$ and
$D^-(f)$ respectively. Recall from \ref{real} that in the real case,
$f$ has two augmentations $g=A_F^+f$ and $\tilde{g}=A_F^-f$ with stable
perturbations $g_\mu(\l,x)=(\l,f_{\lambda^2+\mu}(x))$ and $\tilde
{g}_\mu(\l,x)=(\l,f_{-\lambda^2+\mu}(x))$ respectively.

\begin{prop}

With the above notation
$$\begin{array}{lll}
(i) D^+(g) \cong D^+(f) & &
(ii) D^-(g) \cong S(D^-(f))\\
& & \\
(iii) D^+(\tilde{g}) \cong S(D^+(f)) & & (iv)
D^-(\tilde{g}) \cong
D^-(f)
\end{array} $$
In particular, if $f$ has a good real perturbation then so does
one of its two augmentations.
\end{prop}
\Proof By symmetry it is sufficient to show just the first two homotopy
equivalences. Case (ii) is analogous to Theorem~\ref{suspdis} but
if we
follow the same proof in (i), then since $-\mu_0$ has no real square roots, 
$D^+(g_\mu)$ is a fibre bundle over $B_{\epsilon''}(0)$ with fibre
$D^+(f_\lambda)$. But the total space of a bundle over a contractible
space is
homotopy equivalent to the fibre. \eop

We now describe the topology of a discriminant of a 
stable perturbation over $\RR$ of a real germ in
the normal forms of Theorem~\ref{classi}
(see Remark~\ref{real}). 

Proposition~\ref{tipo} holds in a slightly different
version. Here
we have to consider the two discriminants of a stable perturbation of $h$
as well
as the two discriminants of a stable perturbation
of $f_0$. We leave the straightforward details to the reader,
although we recall that in \ref{goodreal1} we have already shown that
if $f_0$ has a good real perturbation then so does $C_k(f_0)$.
Finally, although we have made no attempt to determine the
number of inequivalent real forms of a binary concatenations of
two real ${\cal A}_e$-codimension
1 multi-germs, the proof of Proposition
\ref{tipoiii} shows
\begin{prop}\label{real3}
Suppose $h=B(f_0,g_0)$ is a binary concatenation of two real multi-germs.
Then the discriminant of a stable perturbation of $h$ (over $\RR$) is
homotopy-equivalent to one of the following four spaces:
$$\begin{array}{ccc}
S(D^+(f_0)\ast D^+(g_0)) & & S(D^-(f_0)\ast D^+(g_0))\\
& & \\
S(D^+(f_0)\ast D^-(g_0)) & & S(D^-(f_0)\ast D^-(g_0)). \end{array} $$
In particular, if $f_0$ and $g_0$ have good real perturbations, then so
does at least one real form of $B(f_0,g_0)$.\eop
\end{prop}

\begin{example} 
{\em
Consider the bi-germ consisting of two prisms on Whitney
cusps, each
transverse to the analytic stratum of the other:
$$h:\cases {f(\lambda, x,\mu)=(\lambda,x^3+\mu x,\mu) \cr
g(z,\delta,\mu)=(z^3-\mu z,\delta,\mu)\cr}$$
The discriminant of each is the product with a line 
of a plane first-order cusp.
The real discriminant of a stable perturbation $h_t$ of $h$
(in which $t$ is added to 
the third component of $g$) is thus the union of two prisms, drawn 
with dotted lines in Figure 5 on page;
its homology is carried by the curvilinear tetrahedron drawn with a 
solid line.  

The intersection of $D(h_t)$ with the horizontal plane
$L_\mu$, for $0 < \mu < t$, is
the union of two pairs of parallel lines, 
$\RR\times D(\tilde{f}_\mu)$ and $D(\tilde{g}_\mu)\times\RR$ (since each of 
$D(\tilde{f}_\mu)$ and $D(\tilde{g}_\mu)$ consists just of a pair of points). 
$L_\mu\cap D(h_\mu)$ retracts to a rectangle,
the intersection of $L_\mu$ with the
(boundary of the) curvilinear tetrahedron. This rectangle is  
the join of $D(\tilde{f}_\mu)$ and $D(\tilde{g}_\mu)$. 
}
\end{example}
\section{Proofs of the main theorems}
\begin{teo}\label{quasi}

Let $h\!:\!(\Bbb {C}^n,T) \to (\Bbb {C}^p,0)$ ($n \geq p-1$, $(n,p)$ nice
dimensions) be a multi-germ of ${\cal {A}}_e$-codimension 1 and corank 1.
Then $h$
is quasihomogeneous.

\end{teo}
\Proof We may suppose $h$ primitive and ignore any submersive branches.  
The proof is by induction on the number, $|T|$, of components of $h$.\\

\ni If $|T|=1$, $h$ is quasihomogeneous by results of Victor Goryunov in
\cite{Goryunov} when $n \geq p$
and by our Proposition~\ref{classification} when $p=n+1$.\\
 
\ni Suppose $h=\{f,g\}$ has more than one branch. If $g$ is not
transverse to $\tau(f)$, then by \ref{classi}, either $f$ and $g$ are both
prisms on Morse singularities, or $h$ is equivalent to $C_k(f_0)$
for some ${\cal A}_e$-codimension 1 germ $f_0$. In the first case $h$ is
plainly quasihomogeneous. In the second, we apply 
the inductive hypothsis to conclude that $f_0$ is quasihomogeneous.
Since we are in the nice dimensions, $f_0$ has a quasihomogeneous
versal unfolding $\tilde{f}$, and by \ref{morpri}, 
$h$ is equivalent to $C_k(f_0)$. Clearly this is quasihomogeneous.

If $f$ is transverse to $\tau(g)$
and vice versa,
then by \ref{banana} the pullback $f_0$ 
of $f$ by $\tau(g)$, and the pull-back $g_0$ of $g$ by $\tau(f)$, 
both have codimension 1. By the induction hypothesis, $f_0$ and $g_0$
are both quasihomogeneous. By
\ref{bin}, $\{f,g\}$
is equivalent to $B(f_0,g_0)$; again, as we are in the nice dimensions, 
$f_0$ and $g_0$ have weighted homogeneous 
${\cal A}_e$-versal unfoldings with unfolding
parameter with positive weight; a representative of $B(f_0,g_0)$
constructed from these ingredients is evidently weighted homogeneous.
\eop
In the next result, we do not distinguish between $\mu_I$ and $\mu_\Delta$,
for the reasons described at the start of Section \ref{concatenation}.

\begin{teo}\label{inumber}
If $h\!:\!(\Bbb {C}^n,T) \to (\Bbb {C}^p,0)$ ($n\geq p-1$, $n,p$ nice 
dimensions)
has corank 1 and ${\cal
{A}}_e$-codimension 1 then $\mu_\Delta(h)=1$ (and in particular $\mu_I=1$
for pure-dimensional multi-germs $\CC^n\to\CC^{n+1}$). 

\end{teo}
\Proof The proof follows exactly the same scheme as the preceding
proof. The starting point for the induction is
now the fact that mono-germs of ${\cal A}_e$-codimension 1 have $\mu_I$
or $\mu_\Delta$ equal to 1, by our proposition \ref{rep} for $n=p-1$,
and by the fact (proved in \cite{damon-mond}) 
that $\mu_\Delta={\cal A}_e$-codimension in the nice dimensions,
for quasi-homogeneous germs $\CC^n\to\CC^p$ with $n\geq p$.

We may suppose $h$ primitive; for by \ref{suspdis}
$D((Ah)_t)\simeq S(D(h_t))$, where the suffix $t$ indicates stable perturbation
and $S$ is suspension.

Since the result is already proven in case all branches have $n\geq p$,
we assume at least one branch has $n=p-1$. Hence by induction and 
\ref{classi} $h$ is equivalent either to $C_k(f_0)$ or to $B(f_0,g_0)$,
where $f_0$ and $g_0$ are quasihomogeneous ${\cal A}_e$-codimension
1 germs. The conclusion now follows by Theorem \ref{conc1}(2) for $C_k(f_0)$
and by Theorem \ref{conc2}(3) for $B(f_0,g_0)$.\eop

\begin{teo}\label{grp}

Let $h\!:\!(\Bbb {C}^n,T) \to (\Bbb {C}^p,0)$ ($n \geq p-1$, $(n,p)$ nice
dimensions) be a multi-germ of ${\cal {A}}_e$-codimension 1 and corank 1. Then
there exists a real form with a good real perturbation.

\end{teo}
\Proof Again, the proof is by induction on $|T|$. 
The result is proven for mono-germs in \cite{mond} (for $n\geq p$) and in
\ref{rep} above for the case $p=n+1$. 
The inductive steps follow, using the classification theorem 
\ref{classi}, by \ref{goodreal1}, \ref{goodreal2} and \ref{real3}.
\eop

{\sc Author addresses}
\begin{tabbing}
\footnotesize{Thomas Cooper}\hskip 1.2in\= 
\footnotesize{David Mond} \hskip 1.2in \= 
\footnotesize{Roberta Wik Atique}\\
\footnotesize{Mathematics Institute}\>\footnotesize{Mathematics Institute}\>
\footnotesize{Instituto 
de Ci\^{e}ncias Matem\^{a}ticas}\\
\footnotesize{University of Warwick}\>
\footnotesize{University of Warwick}\>
\footnotesize{e de Computa\c{c}\~{a}o}\\
\footnotesize{Coventry CV4 7AL}\>
\footnotesize{Coventry CV4 7AL}\>\footnotesize{Caixa Postal 668, 
Sao Carlos, SP}\\
\footnotesize{United Kingdom} \>
\footnotesize{United Kingdom}\>\footnotesize{CEP 13560-970, Brazil}\\\\
\>\footnotesize{mond@maths.warwick.ac.uk}\>\footnotesize{rwik@icmc.sc.usp.br}

\end{tabbing}
\end{document}